\documentclass[12pt,twoside]{article}
\usepackage[hmargin=0.8in,vmargin=0.9in]{geometry}
\usepackage{fancyhdr}
\usepackage{graphicx}
\usepackage{subfigure}
\usepackage{amssymb}
\usepackage{amsmath}
\usepackage{amsfonts}
\usepackage{theorem}
\usepackage{mathrsfs}
\usepackage{mathtools}
\usepackage{bm}
\usepackage{color}
\usepackage{setspace}
\usepackage{exscale}
\usepackage{relsize}
\usepackage{epstopdf}
\usepackage{float}
\usepackage{cite}
\usepackage{nicefrac}
\usepackage{setspace}

\usepackage{booktabs,multirow} 
\usepackage{array} 
\usepackage{paralist} 
\usepackage{verbatim} 
\usepackage{subfigure} 

\theoremstyle{plain}			
\newtheorem{thm}{Theorem}[section]
\newtheorem{lemma}[thm]{Lemma}

{\theorembodyfont{\rmfamily}\newtheorem{remark}{Remark}[section]}

\usepackage{tikz}
\usetikzlibrary{positioning}
\usepackage{xcolor}

\allowdisplaybreaks[1]

\numberwithin{equation}{section}
\numberwithin{figure}{section}
\numberwithin{table}{section}

\newcommand\eref[1]{(\ref{#1})}

\newcommand*\xbar[1]{%
  \hbox{%
    \vbox{%
      \hrule height 0.5pt 
      \kern0.4ex
      \hbox{%
        \kern-0.05em
        \ensuremath{#1}%
        \kern-0.00em
      }%
    }%
  }%
}

\newcommand{\bmF}{\bm{\mathcal{F}}}
\newcommand{\bmG}{\bm{\mathcal{G}}}
\newcommand{\mF}{\bm{F}}

\newcommand{\mG}{\bm{G}}

\newcommand{\mD}{\bm{D}}

\newcommand{\mU}{\bm{U}}

\newcommand{\mo}{\bm{0}}

\newcommand{\dx}{\Delta x}
\newcommand{\dy}{\Delta y}

\newcommand{\hf}{{\frac{1}{2}}}

\newcommand{\jph}{{j+\frac{1}{2}}}
\newcommand{\jmh}{{j-\frac{1}{2}}}
\newcommand{\kph}{{k+\frac{1}{2}}}
\newcommand{\kmh}{{k-\frac{1}{2}}}

\newcommand{\ajphp}{{a_{j+\frac{1}{2}}^+}}
\newcommand{\ajphm}{{a_{j+\frac{1}{2}}^-}}

\def\softd{{\leavevmode\setbox1=\hbox{d}%
          \hbox to 1.05\wd1{d\kern-0.4ex{\char039}\hss}}}

\title{Local Characteristic Decomposition Based Central-Upwind Scheme}
\author{Alina Chertock\thanks{Department of Mathematics, North Carolina State University, Raleigh, NC 27695, USA;
{\tt chertock@math.ncsu.edu}},~Shaoshuai Chu\thanks{Department of Mathematics, Southern University of Science and Technology, Shenzhen,
518055, China; {\tt chuss2019@mail.sustech.edu.cn}},~Michael Herty\thanks{Department of Mathematics, RWTH Aachen University, 52056 Aachen,
Germany; {\tt herty@mathc.rwth-aachen.de}},~Alexander Kurganov\thanks{Department of Mathematics, SUSTech International Center for
Mathematics and Guangdong Provincial Key Laboratory of Computational Science and Material Design, Southern University of Science and
Technology, Shenzhen, 518055, China; {\tt alexander@sustech.edu.cn}},\\~and M\'{a}ria Luk\'{a}\v{c}ov\'{a}-Medvi{\softd}ov\'{a}
\thanks{Institute of Mathematics, University of Mainz, Germany; {\tt lukacova@uni-mainz.de}}}
\date{}
\begin{document}

\maketitle
\begin{abstract}
We propose novel less diffusive schemes for conservative one- and two-dimensional hyperbolic systems of nonlinear partial differential
equations (PDEs). The main challenges in the development of accurate and robust numerical methods for the studied systems come from the
complicated wave structures, such as shocks, rarefactions and  contact discontinuities, arising even for smooth initial conditions. In order
to reduce the diffusion in the original central-upwind schemes, we use a local characteristic decomposition procedure to develop a new
class of central-upwind schemes. We apply the developed schemes to the one- and two-dimensional Euler equations of gas dynamics to
illustrate the performance on a variety of examples. The obtained numerical results clearly demonstrate that the proposed new schemes
outperform the original central-upwind schemes.
\end{abstract}

\noindent
{\bf Key words:} Local characteristic decomposition; central-upwind schemes; hyperbolic systems of conservative laws; Euler equations of gas
dynamics

\noindent
{\bf AMS subject classification:} 76M12, 65M08, 76N15, 35L65, 35L67.

\section{Introduction}
This paper focuses on numerical solutions of hyperbolic systems of conservation laws, which in the two-dimensional (2-D) case, read as
\begin{equation}
\mU_t+\mF(\mU)_x+\mG(\mU)_y=\bm0,
\label{1.1}
\end{equation}
where $x$ and $y$ are spatial variables, $t$ is the time, $\mU\in\mathbb R^d$ is a vector of unknowns, $\mF:\mathbb R^d\to\mathbb R^d$ and
$\mG:\mathbb R^d\to\mathbb R^d$ are the $x$- and $y$-flux functions, respectively.

It is well-known that the solutions of \eref{1.1} may develop complicated wave structures including shocks, rarefactions, and contact
discontinuities even for infinitely smooth initial data, and thus developing highly accurate and robust shock-capturing numerical methods
for solving \eref{1.1} is a challenging task.

Since the pioneering works \cite{Fri,Lax,Godunov59}, a large number of various methods had been introduced; see, e.g., the monographs and
review papers \cite{Hesthaven18,Leveque02,KLR20,Shu20,Tor,BAF} and references therein. Here, we focus on finite-volume (FV) methods, in
which solutions are realized in terms of their cell averages and evolved in time according to the following algorithm. First, a piecewise
polynomial interpolant is reconstructed out of the given cell averages, and then the evolution step is performed using the integral form of
\eref{1.1}. To this end, a proper set of time-space control volumes has to be selected. Depending on this selection, one may distinguish
between two basic classes of finite-volume methods: upwind and central schemes. In upwind schemes, the spatial part of the control volumes
coincides with the FV cells and therefore, one needs to (approximately) solve (generalized) Riemann problems at every cell interface; see,
e.g., \cite{Godunov59,Tor,BAF} and references therein. This helps the upwind schemes to achieve very high resolution. At the same time, it
might be quite complicated and even impossible to solve the (generalized) Riemann problems for general systems of conservation laws. Central
schemes offer an attractive simple alternative to the upwind ones. In staggered central schemes, first proposed for one-dimensional (1-D)
systems in \cite{Nessyahu90} and then extended to higher order \cite{Levy99,Liu98} and multiple number of dimensions
\cite{Lie03a,Jiang98,Arminjon95}, the control volumes use the staggered spatial parts so that the FV cell interfaces remain inside the
control volumes. This allows to avoid even approximately solving any Riemann problems, which makes staggered central schemes easy to
implement for a wide variety of hyperbolic systems. The major drawback of staggered central schemes is, however, their relatively large
numerical dissipation as they basically average over the Riemann fans rather than resolving them.

In order to reduce the amount of excessive numerical dissipation present in central schemes, a class of central-upwind (CU) schemes have
been proposed in \cite{Kurganov01,Kurganov00}. These schemes are based on nonuniform control volumes, whose spatial size taken to be
proportional to the local speeds of propagation. This allows one to minimize the area over which the solution averaging occurs still
without (approximately) solving any (generalized) Riemann problems. The upwind features of the CU schemes can be seen, for example, when
they are applied to simpler systems. For instance, the CU scheme from \cite{Kurganov01} reduces to the upwind one when it is applied to a
system whose Jacobian contains only positive (only negative) eigenvalues. Another advantage of the CU schemes is related to the fact that
unlike the staggered central schemes, they admit a particularly simple semi-discrete form. This observation is the basis of the
modifications of the CU schemes we propose in this paper.

Even though the CU schemes from \cite{Kurganov01,Kurganov00} are quite accurate, efficient and robust tools for a wide variety of hyperbolic
systems, higher resolution of the numerical solutions can be achieved by further reducing numerical dissipation. This can be done in a
number of different ways, for example: (i) by introducing a more accurate evolution procedure, which leads to a ``built-in'' anti-diffusion
term \cite{Kurganov07}; (ii) by implementing a numerical dissipation switch to control the amount of numerical dissipation present in the CU
schemes \cite{Kurganov21a}; (iii) by obtaining more accurate estimates for the one-sided local speeds of propagation using the discrete
Rankine-Hugoniot conditions \cite{Garg21}.

Another way to control the amount of numerical dissipation present in the CU schemes or any other FV methods is by adjusting the nonlinear
limiting mechanism used in the piecewise polynomial reconstruction. It is well-known that sharper reconstructions may lead to larger
numerical oscillations and a way to reduce these oscillations is to reconstruct the characteristic variables rather than the conservative
ones; see, e.g., \cite{Qiu02}. This can be done using the local characteristic decomposition (LCD); see, e.g.,
\cite{don9,JSZ,Liu17,Nonomura20,Qiu02,Shu20,wang18} and references therein.

In this paper, we modify the CU schemes from \cite{Kurganov01} by applying  LCD to the numerical diffusion part of the schemes. The
obtained new LCD-based CU schemes contain substantially smaller amount of numerical dissipation, which leads to a significantly improved
resolution of the computed solution compared with the original CU schemes. As observed above, the key idea is that the new LCD-based CU
scheme reduces to the upwind scheme when applied to 1-D linear hyperbolic systems
\begin{equation}
\mU_t+A\mU_x=\mo,
\label{1.2}
\end{equation}
where $A$ is a constant matrix (disregarding the sign of the eigenvalues of $A$). This feature suggests that the new CU schemes have more
built-in upwinding compared with their predecessors.

The paper is organized as follows. In \S\ref{sec2}, we briefly describe the 1-D second-order FV CU scheme from \cite{Kurganov01}. In
\S\ref{sec3}, we introduce the proposed new 1-D  LCD-based CU scheme and show that the developed 1-D  scheme reduces to the upwind scheme
when applied to the linear system \eref{1.2}. In \S \ref{sec5}, we construct the 2-D LCD-based CU scheme. Finally, in \S\ref{sec6}, we test
the proposed schemes on a number of 1-D and 2-D numerical examples for the Euler equations of gas dynamics. We demonstrate high accuracy,
efficiency, stability, and robustness of the new LCD-based CU schemes, which clearly outperform the second-order CU scheme from
\cite{Kurganov01,Kurganov02}.

\section{1-D Central-Upwind Scheme: A Brief Overview}\label{sec2}
In this section, we consider the 1-D  hyperbolic system of conservation laws
\begin{equation}
\mU_t+\mF(\mU)_x=\bm0,
\label{2.1}
\end{equation}
and describe the second-order semi-discrete CU scheme from \cite{Kurganov01}. To this end, we assume that the computational domain is
covered with the uniform cells $C_j:=[x_\jmh,x_\jph]$ of size $\dx$ centered at $x_j=(x_\jmh+x_\jph)/2$ and denote by $\xbar\mU_j(t)$ cell
averages of $\mU(\cdot,t)$ over the corresponding intervals $C_j$, that is,
\begin{equation*}
\xbar\mU_j(t):\approx\frac{1}{\dx}\int\limits_{C_j}\mU(x,t)\,{\rm d}x.
\end{equation*}
We also assume that at certain time $t\ge0$, the cell average values $\xbar\mU_j$ are available and from here on we suppress the
time-dependence of all of the indexed quantities for the sake of brevity.

According to \cite{Kurganov01}, the computed cell averages \eref{2.1} are evolved in time by solving the following system of ordinary
differential equations (ODEs):
\begin{equation}
\frac{{\rm d}\xbar\mU_j}{{\rm d}t}=-\frac{\bmF_\jph-\bmF_\jmh}{\dx},
\label{2.2}
\end{equation}
where $\bmF_\jph$ are the CU numerical fluxes given by
\begin{equation}
\bmF_\jph=\frac{\ajphp\mF^-_\jph-\ajphm\mF^+_\jph}{\ajphp-\ajphm}+\frac{\ajphp\ajphm}{\ajphp-\ajphm}\left(\mU^+_\jph-\mU^-_\jph\right).
\label{2.3}
\end{equation}
Here, $\mF^\pm_\jph:=\mF\big(\mU^\pm_\jph\big)$ and $\mU^\pm_\jph$ are the right/left-sided values reconstructed out of the given set of
cell averages $\{\xbar\mU_j\}$. The one-sided local speeds  of propagation $a^\pm_\jph$ are estimated using the largest and the smallest
eigenvalues of the Jacobian $A(\mU):=\frac{\partial\mF}{\partial\mU}(\mU)$, which we denote by
$\lambda_1(A(\mU))\le\ldots\le\lambda_d(A(\mU))$. This can be done, for example, by taking
\begin{equation}
\begin{aligned}
&a^+_\jph=\max\left\{\lambda_d\big(A(\mU^+_\jph)\big),\lambda_d\big(A(\mU^-_\jph)\big),0\right\},\\
&a^-_\jph=\min\left\{\lambda_1\big(A(\mU^+_\jph)\big),\lambda_1\big(A(\mU^-_\jph)\big),0\right\}.
\end{aligned}
\label{2.4}
\end{equation}

The (formal) order of the scheme \eref{2.2}--\eref{2.4} is determined by the accuracy the point values $\mU^\pm_\jph$ are reconstructed with
and the order of the ODE solver used to integrate the ODE system \eref{2.2} in time. In this paper, we use the LCD-based second-order
piecewise linear reconstruction described in \S\ref{sec2.1} and the three-stage third-order strong stability preserving (SSP) Runge-Kutta
solver; see, e.g., \cite{Gottlieb11,Gottlieb12}.

\subsection{ LCD-Based Piecewise Linear Reconstruction}\label{sec2.1}
In order to ensure at least second-order of accuracy of the CU scheme \eref{2.2}--\eref{2.4}, one has to use a second-order piecewise
linear reconstruction, which will be non-oscillatory provided the numerical derivatives are computed using a nonlinear limiter. A library of
different limiters are available; see, e.g., \cite{Lie03,Nessyahu90,Sweby84,BAF,Hesthaven18,Leveque02,Tor}. The limiters can be classified
as dissipative, compressive or overcompressive; see \cite{Lie03}. When applied to the conservative variables $\mU$, dissipative limiters may
introduce excessive numerical dissipation, which may results in oversmeared solution discontinuities. On the other hand, compressive and
overcompressive limiters tend to lead to quite large oscillations near discontinuities and artificial sharpening smooth parts of the
solution (compressive limiters may lead to kinks, that is, jumps in the first derivatives, while the overcompressive limiters may lead to
jumps in the solution itself); see, e.g., \cite{Lie03}. Part of these difficulties can be overcome by reconstructing local characteristic
variables as it was demonstrated in \cite{Qiu02} in the context of higher-order WENO reconstructions.

In this paper, we use the minmod limiter, which is, according to \cite{Lie03}, the most compressive out of dissipative limiters. As this
limiter typically leads to quite large oscillations when applied to the conservative variables, we implement it in the LCD framework.
Specifically, we first introduce the matrix $\widehat A_\jph=A(\widehat\mU_\jph)$, where $A(\mU)=\frac{\partial\mF}{\partial\mU}(\mU)$ is
the Jacobian and $\widehat\mU_\jph$ is either a simple average $(\xbar\mU_j+\xbar\mU_{j+1})/2$ or another type of average of the
$\xbar\mU_j$ and $\xbar\mU_{j+1}$ states.

As long as the system \eref{1.1} is strictly hyperbolic, we compute the matrices $R_\jph$ and $R^{-1}_\jph$ such that
$R^{-1}_\jph\widehat A_\jph R_\jph$ is a diagonal matrix and introduce the local characteristic variables $\bm\Gamma_k$ in the neighborhood
of $x=x_\jph$:
$$
\bm\Gamma_k=R^{-1}_\jph\xbar\mU_k,\quad k=j-1,\ldots,j+2.
$$

Equipped with the values $\bm\Gamma_{j-1}$, $\bm\Gamma_j$, $\bm\Gamma_{j+1}$ and $\bm\Gamma_{j+2}$, we reconstruct $\bm\Gamma$ in the cell
$C_j$ by computing
\begin{equation}
(\bm\Gamma_x)_j={\rm minmod}\left(2\,\frac{\bm\Gamma_j-\bm\Gamma_{j-1}}{\dx},\,\frac{\bm\Gamma_{j+1}-\bm\Gamma_{j-1}}{2\dx},\,
2\,\frac{\bm\Gamma_{j+1}-\bm\Gamma_j}{\dx}\right),
\label{2.5}
\end{equation}
and
\begin{equation}
(\bm\Gamma_x)_{j+1}={\rm minmod}\left(2\,\frac{\bm\Gamma_{j+1}-\bm\Gamma_{j}}{\dx},\,\frac{\bm\Gamma_{j+2}-\bm\Gamma_{j}}{2\dx},\,
2\,\frac{\bm\Gamma_{j+2}-\bm\Gamma_{j+1}}{\dx}\right),
\label{2.6}
\end{equation}
where the minmod function, defined as
\begin{equation}
{\rm minmod}(z_1,z_2,\ldots):=\begin{cases}
\min_j\{z_j\}&\mbox{if}~z_j>0\quad\forall\,j,\\
\max_j\{z_j\}&\mbox{if}~z_j<0\quad\forall\,j,\\
0            &\text{otherwise,}
\end{cases}
\label{2.7}
\end{equation}
is applied in the component-wise manner. The slopes \eref{2.5} and \eref{2.6} allow to evaluate
$$
\bm\Gamma^-_\jph=\bm\Gamma_j+\frac{\dx}{2}(\bm\Gamma_x)_j\quad\mbox{and}\quad
\bm\Gamma^+_\jph=\bm\Gamma_{j+1}-\frac{\dx}{2}(\bm\Gamma_x)_{j+1},
$$
and  obtain the corresponding point values of $\mU$ by
\begin{equation}
\mU^\pm_\jph=R_\jph\bm\Gamma^\pm_\jph.
\label{2.8}
\end{equation}
\begin{remark}
In Appendix \ref{appa}, we provide a detailed explanation on how the average matrix $\widehat A_\jph$ and the corresponding matrices
$R_\jph$ and $R^{-1}_\jph$ are computed in the case of the Euler equation of gas dynamics.
\end{remark}

\section{1-D LCD-Based Central-Upwind Scheme}\label{sec3}
In this section, we introduce a new 1-D LCD-based CU scheme, in which the amount of numerical dissipation is substantially reduced compared
with the original CU scheme \eref{2.2}--\eref{2.3}. To this end, we first rewrite the numerical flux \eref{2.3} of the original CU scheme in
the following form:
\begin{equation}
\bmF_\jph=\frac{\mF_j+\mF_{j+1}}{2}+\mD_\jph,
\label{3.1}
\end{equation}
where $\mF_j:=\mF(\,\xbar\mU_j)$ and $\mD_\jph$ is the numerical diffusion given by
\begin{equation*}
\begin{aligned}
\mD_\jph&=\frac{a^+_\jph}{a^+_\jph-a^-_\jph}\left[\mF^-_\jph-\frac{\mF_j+\mF_{j+1}}{2}\right]-\frac{a^-_\jph}{a^+_\jph-a^-_\jph}
\left[\mF^+_\jph-\frac{\mF_j+\mF_{j+1}}{2}\right]\\
&+\frac{a^+_\jph a^-_\jph}{a^+_\jph-a^-_\jph}\left(\mU^+_\jph-\mU^-_\jph\right),
\end{aligned}
\end{equation*}
which, in turn, can be rewritten with the help of the matrix $R_\jph$ introduced in \S\ref{sec2.1} and using \eref{2.8} as follows:
\begin{equation}
\begin{aligned}
\mD_\jph&=R_\jph P_\jph R^{-1}_\jph \left[\mF^-_\jph-\frac{\mF_j+\mF_{j+1}}{2}\right]+R_\jph M_\jph R^{-1}_\jph
\left[\mF^+_\jph-\frac{\mF_j+\mF_{j+1}}{2}\right]\\
&+R_\jph Q_\jph\left(\bm\Gamma^+_\jph-\bm\Gamma^-_\jph\right).
\end{aligned}
\label{3.2}
\end{equation}
Here, $P_\jph$, $M_\jph$, and $Q_\jph$ are the diagonal matrices
\begin{equation}
\begin{aligned}
&P_\jph={\rm diag}\big((P_1)_\jph,\ldots,(P_d)_\jph\big),\quad M_\jph={\rm diag}\big((M_1)_\jph,\ldots,(M_d)_\jph\big),\\
&Q_\jph={\rm diag}\big((Q_1)_\jph,\ldots,(Q_d)_\jph\big)
\end{aligned}
\label{3.3}
\end{equation}
with
\begin{equation}
\big((P_i)_\jph,(M_i)_\jph,(Q_i)_\jph\big)=\frac{1}{a^+_{\jph}-a^-_\jph}\big(a^+_\jph,-a^-_\jph,a^+_\jph a^-_\jph\big),\quad i=1,\ldots,d.
\label{3.4}
\end{equation}

The main idea of the new the 1-D LCD-based CU scheme is to replace the constant entries \eref{3.4} in the diagonal matrices $P_\jph$,
$M_\jph$, and $Q_\jph$ with the corresponding characteristic entries of the diagonal matrix $R^{-1}_\jph\widehat A_\jph R_\jph$. This is
motivated by the fact that in the linear case with $\mF(\mU)=A\mU$, the local characteristic speed might be different for each
(diagonalized) component $\bm\Gamma$. Hence, a diffusion that depends on the corresponding local speed of the each component may lead to a
sharper resolution of possible discontinuities.

Proposed modifications lead to the semi-discrete scheme
\begin{equation}
\frac{{\rm d}\xbar\mU_j}{{\rm d}t}=-\frac{\bmF^{\rm LCD}_\jph-\bmF^{\rm LCD}_\jmh}{\dx},
\label{3.5}
\end{equation}
where the numerical fluxes $\bmF^{\rm LCD}_\jph$ are the LCD modifications of \eref{3.1}:
\begin{equation}
\bmF^{\rm LCD}_\jph=\frac{\mF_j+\mF_{j+1}}{2}+\mD^{\rm LCD}_\jph,
\label{3.6}
\end{equation}
with the following LCD modification of the numerical diffusion term \eref{3.2}--\eref{3.4}:
\begin{equation}
\begin{aligned}
\mD^{\rm LCD}_\jph&=R_\jph P^{\rm LCD}_\jph R^{-1}_\jph \left[\mF^-_\jph-\frac{\mF_j+\mF_{j+1}}{2}\right]+
R_\jph M^{\rm LCD}_\jph R^{-1}_\jph\left[\mF^+_\jph-\frac{\mF_j+\mF_{j+1}}{2}\right]\\
&+R_\jph Q^{\rm LCD}_\jph\left(\bm\Gamma^+_\jph-\bm\Gamma^-_\jph\right).
\end{aligned}
\label{3.7}
\end{equation}
Here,
\begin{equation*}
\begin{aligned}
&P^{\rm LCD}_\jph={\rm diag}\big((P^{\rm LCD}_1)_\jph,\ldots,(P^{\rm LCD}_d)_\jph\big),\quad
M^{\rm LCD}_\jph={\rm diag}\big((M^{\rm LCD}_1)_\jph,\ldots,(M^{\rm LCD}_d)_\jph\big),\\
&Q^{\rm LCD}_\jph={\rm diag}\big((Q^{\rm LCD}_1)_\jph,\ldots,(Q^{\rm LCD}_d)_\jph\big)
\end{aligned}
\end{equation*}
with
\begin{align}
&\hspace*{-0.6cm}\big((P^{\rm LCD}_i)_\jph,(M^{\rm LCD}_i)_\jph,(Q^{\rm LCD}_i)_\jph\big)\label{3.8}\\
&=\left\{\begin{aligned}
&\frac{1}{(\lambda^+_i)_\jph-(\lambda^-_i)_\jph}\big((\lambda^+_i)_\jph,-(\lambda^-_i)_\jph,(\lambda^+_i)_\jph(\lambda^-_i)_\jph\big)&&
\mbox{if}~(\lambda^+_i)_\jph-(\lambda^-_i)_\jph> \varepsilon,\\
&0&&\mbox{otherwise},
\end{aligned}\right.\nonumber
\end{align}
and
\begin{equation*}
\begin{aligned}
(\lambda^+_i)_\jph&=\max\left\{\lambda_i\big(A(\mU^-_\jph)\big),\,\lambda_i\big(A(\mU^+_\jph)\big),\, 0\right\},\\
(\lambda^-_i)_\jph&=\min\left\{\lambda_i\big(A(\mU^-_\jph)\big),\,\lambda_i\big(A(\mU^+_\jph)\big),\, 0\right\},
\end{aligned}
\end{equation*}
for $i=1,\ldots,d$.

Finally, $\varepsilon$ in \eref{3.8} is a very small desingularization constant, taken $\varepsilon=10^{-10}$ in all of the numerical
examples reported in \S\ref{sec6}.

It is crucial to note that the proposed 1-D LCD-based CU scheme reduces to the second-order semi-discrete upwind scheme when applied to a
linear hyperbolic system \eref{1.2} with constant coefficients. This is proven in the following lemma.
\begin{lemma}
If
\begin{equation}
\mF(\mU)=A\mU,
\label{3.9}
\end{equation}
where $A$ is a constant matrix, then the scheme \eref{3.5}--\eref{3.7} becomes the second-order semi-discrete upwind scheme.
\end{lemma}
{\bf Proof.} We note that in the linear case with constant coefficients, the corresponding matrix composed of the right eigenvectors of $A$
is also a constant matrix, that is, $R_\jph\equiv R$, and matrices $P^{\rm LCD}_\jph$ and $M^{\rm LCD}_\jph$ reduce to
\begin{equation}
\begin{aligned}
&P^{\rm LCD}_\jph\equiv P={\rm diag}\left(\max\{{\rm sign}(\lambda_1),0\},\ldots,\max\{{\rm sign}(\lambda_d),0\}\right),\\
&M^{\rm LCD}_\jph\equiv M={\rm diag}\left(\min\{{\rm sign}(\lambda_1),0\},\ldots,\min\{{\rm sign}(\lambda_d),0\}\right),
\end{aligned}
\label{3.10}
\end{equation}
while the matrix $Q^{\rm LCD}_\jph\equiv0$ as $\lambda^+_i\lambda^-_i=\max(\lambda_i,0)\cdot\min(\lambda_i,0)=0$.

Taking into account the above simplifications and substituting \eref{3.9}, \eref{3.10} into \eref{3.7}, yields the following expression for
the numerical diffusion term:
\begin{equation}
\begin{aligned}
&\mD^{\rm LCD}_\jph=RPR^{-1}A\left[\mU^-_\jph-\frac{\xbar\mU_j+\xbar\mU_{j+1}}{2}\right]
-RM R^{-1}A\left[\mU^+_\jph-\frac{\xbar\mU_j+\xbar\mU_{j+1}}{2}\right]\\
&=RPR^{-1}AR\left[R^{-1}\left(\mU^-_\jph-\frac{\xbar\mU_j+\xbar\mU_{j+1}}{2}\right)\right]
-RMR^{-1}AR\left[R^{-1} \left(\mU^+_\jph-\frac{\xbar\mU_j+\xbar\mU_{j+1}}{2}\right)\right]\\
&=R\Lambda^+\left(\bm\Gamma^-_\jph-\frac{\bm\Gamma_j+\bm\Gamma_{j+1}}{2}\right)
+R\Lambda^-\left(\bm\Gamma^+_\jph-\frac{\bm\Gamma_j+\bm\Gamma_{j+1}}{2}\right),
\end{aligned}
\label{3.11}
\end{equation}
where
\begin{equation*}
\begin{aligned}
\Lambda^+&=PR^{-1}AR=P\Lambda={\rm diag}\left(\max\{\lambda_1,0\},\ldots,\max\{\lambda_d,0\}\right),\\
\Lambda^-&=MR^{-1}AR=M\Lambda={\rm diag}\left(\min\{\lambda_1,0\},\ldots,\min\{\lambda_d,0\}\right).
\end{aligned}
\end{equation*}
Finally, combining \eref{3.9} and \eref{3.11} together and including them into \eref{3.6}, we obtain
\begin{equation}
\begin{aligned}
\bmF^{\rm LCD}_\jph&=\hf A\left(\,\xbar\mU_j+\xbar\mU_{j+1}\right)+R\Lambda^+\left(\bm\Gamma^-_\jph-
\frac{\bm\Gamma_j+\bm\Gamma_{j+1}}{2}\right)+R\Lambda^-\left(\bm\Gamma^+_\jph-\frac{\bm\Gamma_j+\bm\Gamma_{j+1}}{2}\right)\\
&\hspace*{-0.6cm}=\hf(A^++A^-)\left(\,\xbar\mU_j+\xbar\mU_{j+1}\right)+A^+\left(\mU^-_\jph-\frac{\mU_j+\mU_{j+1}}{2}\right)+
A^-\left(\mU^+_\jph-\frac{\mU_j+\mU_{j+1}}{2}\right)\\
&\hspace*{-0.6cm}=A^+\mU^-_\jph+A^-\mU^+_\jph,
\end{aligned}
\label{3.12}
\end{equation}
where $A^\pm=R\Lambda^\pm R^{-1}$. This confirms that the scheme \eref{3.5}, \eref{3.12} is the second-order semi-discrete upwind scheme.
$\hfill\blacksquare$

\section{2-D LCD-Based Central-Upwind Scheme}\label{sec5}
In this section, we generalize the 1-D LCD-based CU scheme introduced in \S\ref{sec3} for the 2-D hyperbolic system of conservation laws
\eref{1.1}. We design the 2-D LCD-based CU scheme in a ``dimension-by-dimension'' manner, so that it reads as
\begin{equation*}
\frac{{\rm d}\xbar\mU_{j,k}}{{\rm d}t}=-\frac{\bmF^{\rm LCD}_{\jph,k}-\bmF^{\rm LCD}_{\jmh,k}}{\dx}-
\frac{\bmG^{\rm LCD}_{j,\kph}-\bmG^{\rm LCD}_{j,\kmh}}{\dy},
\end{equation*}
where
\begin{equation*}
\bmF^{\rm LCD}_{\jph,k}=\frac{\mF_{j,k}+\mF_{j+1,k}}{2}+\mD^{\rm LCD}_{\jph,k},\quad\bmG^{\rm LCD}_{j,\kph}=
\frac{\mG_{j,k}+\mG_{j,k+1}}{2}+\mD^{\rm LCD}_{j,\kph}.
\end{equation*}
Here, $\mF_{j,k}:=\mF(\,\xbar\mU_{j,k})$, $\mG_{j,k}:=\mG(\,\xbar\mU_{j,k})$, and $\mD^{\rm LCD}_{\jph,k}$ and $\mD^{\rm LCD}_{j,\kph}$ are
the numerical diffusion terms defined by
\begin{equation*}
\begin{aligned}
\mD^{\rm LCD}_{\jph,k}&=R_{\jph,k}P^{\rm LCD}_{\jph,k}R^{-1}_{\jph,k}\left[\mF^{\rm E}_{j,k}-\frac{\mF_{j,k}+\mF_{j+1,k}}{2}\right]\\
&+R_{\jph,k}M^{\rm LCD}_{\jph,k}R^{-1}_{\jph,k}\left[\mF^{\rm W}_{j+1,k}-\frac{\mF_{j,k}+\mF_{j+1,k}}{2}\right]+
R_{\jph,k}Q^{\rm LCD}_{\jph,k}\left(\bm\Gamma^{\rm W}_{j+1,k}-\bm\Gamma^{\rm E}_{j,k}\right),\\
\mD^{\rm LCD}_{j,\kph}&=R_{j,\kph}P^{\rm LCD}_{j,\kph}R^{-1}_{j,\kph}\left[\mG^{\rm N}_{j,k}-\frac{\mG_{j,k}+\mG_{j,k+1}}{2}\right]\\
&-R_{j,\kph}M^{\rm LCD}_{j,\kph}R^{-1}_{j,\kph}\left[\mG^{\rm S}_{j,k+1}-\frac{\mG_{j,k}+\mG_{j,k+1}}{2}\right]+
R_{j,\kph}Q^{\rm LCD}_{j,\kph}\left[\bm\Gamma^{\rm S}_{j,k+1}-\bm\Gamma^{\rm N}_{j,k}\right].
\end{aligned}
\end{equation*}
The matrices $R_{\jph,k}$, $R^{-1}_{\jph,k}$ and $R_{j,\kph}$, $R^{-1}_{j,\kph}$ are the matrices such that
$R^{-1}_{\jph,k}\widehat A_{\jph,k}R_{\jph,k}$ and $R^{-1}_{j,\kph}\widehat B_{j,\kph}R_{j,\kph}$ are diagonal. Here,
$\widehat A_{\jph,k}=A(\widehat\mU_{\jph,k})$, $\widehat B_{j,\kph}=B(\widehat\mU_{j,\kph})$ with
$A(\mU)=\frac{\partial\mF(\mU)}{\partial\mU}$, $B(\mU)=\frac{\partial\mG(\mU)}{\partial\mU}$, and $\widehat\mU_{\jph,k}$,
$\widehat\mU_{j,\kph}$ are either simple averages $(\,\xbar\mU_{j,k}+\xbar\mU_{j+1,k})/2$, $(\,\xbar\mU_{j,k}+\xbar\mU_{j,k+1})/2$ or
another type of averages of $\xbar\mU_{j,k}$, $\xbar\mU_{j+1,k}$ and $\xbar\mU_{j,k}$, $\xbar\mU_{j,k+1}$ states, respectively. The
numerical fluxes $\mF^{\rm E,W}_{j,k}$ and $\mG^{\rm N,S}_{j,k}$ are defined by $\mF^{\rm E,W}_{j,k}:=\mF(\mU^{\rm E,W}_{j,k})$,
$\mG^{\rm N,S}_{j,k}:=\mG(\mU^{\rm N,S}_{j,k})$, and the details of reconstructing the point values $\mU^{\rm E,W}_{j,k}$ and
$\mU^{\rm N,S}_{j,k}$, $\bm\Gamma^{\rm E,W}_{j,k}$ and $\bm\Gamma^{\rm N,S}_{j,k}$ for the 2-D Euler equations of gas dynamics are provided
in Appendix \ref{appb}. Finally, the diagonal matrices $P^{\rm LCD}_{\jph,k}$, $M^{\rm LCD}_{\jph,k}$, $Q^{\rm LCD}_{\jph,k}$,
$P^{\rm LCD}_{j,\kph}$, $M^{\rm LCD}_{j,\kph}$, and $Q^{\rm LCD}_{j,\kph}$ are defined by
\begin{equation*}
\begin{aligned}
&P^{\rm LCD}_{\jph,k}={\rm diag}\left(\big(P^{\rm LCD}_1\big)_{\jph,k},\ldots,\big(P^{\rm LCD}_d\big)_{\jph,k}\right),&&
P^{\rm LCD}_{j,\kph}={\rm diag}\left(\big(P^{\rm LCD}_1\big)_{j,\kph},\ldots,\big(P^{\rm LCD}_d\big)_{j,\kph}\right),\\
&M^{\rm LCD}_{\jph,k}={\rm diag}\left(\big(M^{\rm LCD}_1\big)_{\jph,k},\ldots,\big(M^{\rm LCD}_d\big)_{\jph,k}\right),&&
M^{\rm LCD}_{j,\kph}={\rm diag}\left(\big(M^{\rm LCD}_1\big)_{j,\kph},\ldots,\big(M^{\rm LCD}_d\big)_{j,\kph}\right),\\
&Q^{\rm LCD}_{\jph,k}={\rm diag}\left(\big(Q^{\rm LCD}_1\big)_{\jph,k},\ldots,\big(Q^{\rm LCD}_d\big)_{\jph,k}\right),&&
Q^{\rm LCD}_{j,\kph}={\rm diag}\left(\big(Q^{\rm LCD}_1\big)_{j,\kph},\ldots,\big(Q^{\rm LCD}_d\big)_{j,\kph}\right),
\end{aligned}
\end{equation*}
where for all $i=1,\ldots,d$
$$
\begin{aligned}
&\hspace*{-1.0cm}\left((P^{\rm LCD}_i)_{\jph,k},(M^{\rm LCD}_i)_{\jph,k},(Q^{\rm LCD}_i)_{\jph,k}\right)\\
&=\left\{\begin{aligned}
&\frac{1}{\Delta(\lambda_i)_{\jph,k}}
\left((\lambda^+_i)_{\jph,k},-(\lambda^-_i)_{\jph,k},(\lambda^+_i)_{\jph,k}(\lambda^-_i)_{\jph,k}\right)
&&\mbox{if}~\Delta(\lambda_i)_{\jph,k}>\varepsilon,\\
&0&&\mbox{otherwise},
\end{aligned}\right.\\
&\hspace*{-1.0cm}\left((P^{\rm LCD}_i)_{j,\kph},(M^{\rm LCD}_i)_{j,\kph},(Q^{\rm LCD}_i)_{j,\kph}\right)\\
&=\left\{\begin{aligned}
&\frac{1}{\Delta(\mu_i)_{j,\kph}}\left((\mu^+_i)_{j,\kph},-(\mu^-_i)_{j,\kph},(\mu^+_i)_{j,\kph}(\mu^-_i)_{j,\kph}\right)
&&\mbox{if}~\Delta(\mu_i)_{j,\kph}>\varepsilon,\\
&0&&\mbox{otherwise}.
\end{aligned}\right.
\end{aligned}
$$
Here, $\Delta(\lambda_i)_{\jph,k}:=(\lambda^+_i)_{\jph,k}-(\lambda^-_i)_{\jph,k}$,
$\Delta(\mu_i)_{j,\kph}:=(\mu^+_i)_{j,\kph}-(\mu^-_i)_{j,\kph}$, and
\begin{equation*}
\begin{aligned}
(\lambda^+_i)_{\jph,k}=\max\left\{\lambda_i\big(A(\mU^{\rm E}_{j,k})\big),\,\lambda_i\big(A(\mU^{\rm W}_{j+1,k})\big),\,0\right\},\\
(\lambda^-_i)_{\jph,k}=\min\left\{\lambda_i\big(A(\mU^{\rm E}_{j,k})\big),\,\lambda_i\big(A(\mU^{\rm W}_{j+1,k})\big),\,0\right\},\\
(\mu^+_i)_{j,\kph}=\max\left\{\mu_i\big(B(\mU^{\rm N}_{j,k})\big),\,\mu_i\big(B(\mU^{\rm S}_{j,k+1})\big),\,0\right\},\\
(\mu^-_i)_{j,\kph}=\min\left\{\mu_i\big(B(\mU^{\rm N}_{j,k})\big),\,\mu_i\big(B(\mU^{\rm S}_{j,k+1})\big),\,0\right\},
\end{aligned}
\end{equation*}
where $\lambda_i$ and $\mu_i$ are the eigenvalues of the Jacobians $A(\mU)$ and $B(\mU)$: $\lambda_1(A)\le\ldots\le\lambda_d(A)$ and
$\mu_1(B)\le\ldots\le\mu_d(B)$, respectively.

\section{Numerical Examples}\label{sec6}
In this section, we apply the proposed LCD-based CU schemes, which will be referred to as the New CU schemes, to the 1-D and 2-D Euler
equations of gas dynamics described in Appendices \ref{appa} and \ref{appb}, respectively. We conduct several numerical experiments and
compare the performance of the {\em New CU} schemes with that of the corresponding 1-D and 2-D second-order CU schemes from
\cite{Kurganov01} and \cite{Kurganov02}, respectively, which will be referred to as the {\em Old CU} schemes.

In Examples 1--8, we take the specific heat ratio $\gamma=1.4$, while in Example 9, we set $\gamma=5/3$. In all of the examples, we use the
CFL number 0.4.

\subsection{One-Dimensional Examples}
\subsubsection*{Example 1---``Shock-Bubble'' Interaction Problem}
\begin{figure}[ht!]
\centerline{
\includegraphics[trim=1.1cm 0.4cm 1.3cm 0.4cm, clip, width=6.4cm]{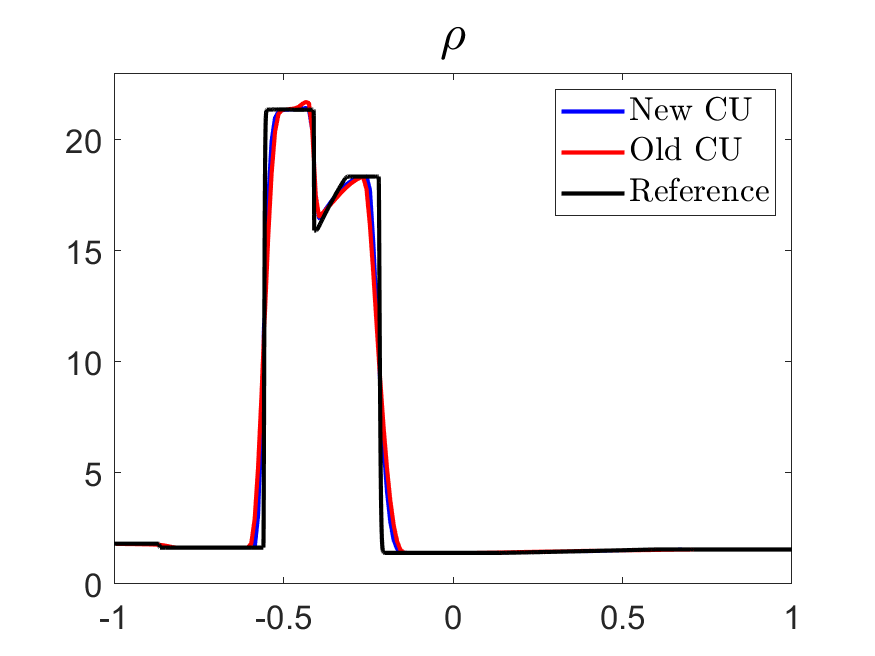}\hspace*{1.0cm}
\includegraphics[trim=1.1cm 0.4cm 1.3cm 0.4cm, clip, width=6.4cm]{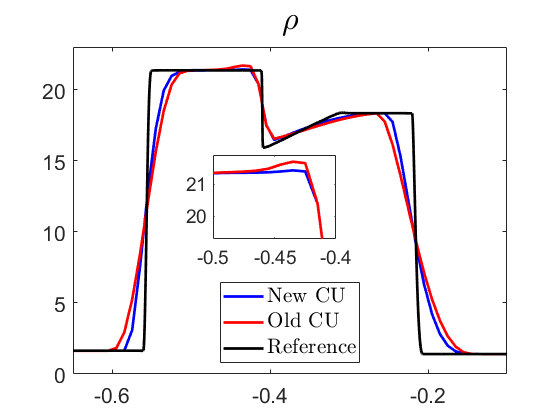}}
\caption{\sf Example 1: Density ($\rho$) computed by the New and Old CU schemes and zoom at $x\in[-0.65,-0.1]$ and $x\in[-0.5,-0.4]$.
\label{fig1}}
\end{figure}

In the first example, we consider the ``shock-bubble'' interaction problem taken from \cite{Kurganov22}. The initial data, given by
\begin{equation*}
(\rho, u,p)(x,0)=\begin{cases}
(13.1538,0,1)&\mbox{if}~|x|<0.25,\\
(1.3333,-0.3535,1.5)&\mbox{if}~x>0.75,\\
(1,0,1)&\mbox{otherwise},
\end{cases}
\end{equation*}
correspond to a left-moving shock, initially located at $x=0.75$, and a bubble of radius 0.25, initially located at the origin. We implement
solid wall boundary conditions on the left and free boundary conditions on the right of the computational domain $[-1,1]$. We compute the
numerical solution until the final time $t=3$ using both he New and Old CU schemes on the uniform mesh with $\dx=1/100$. The obtained
numerical results are presented in Figure \ref{fig1} along with the reference solution computed by the Old CU scheme on a much finer mesh
with $\dx=1/2000$. As one can see, the results obtained by the New CU scheme are sharper and less oscillatory compared to the Old CU scheme
results, which have lower resolution of the shock waves and contain an oscillation at the contact discontinuity.

\subsubsection*{Example 2---Shock-Entropy Wave Interaction Problem}
\begin{figure}[ht!]
\centerline{
\includegraphics[trim=0.9cm 0.4cm 1.3cm 0.4cm, clip, width=6.4cm]{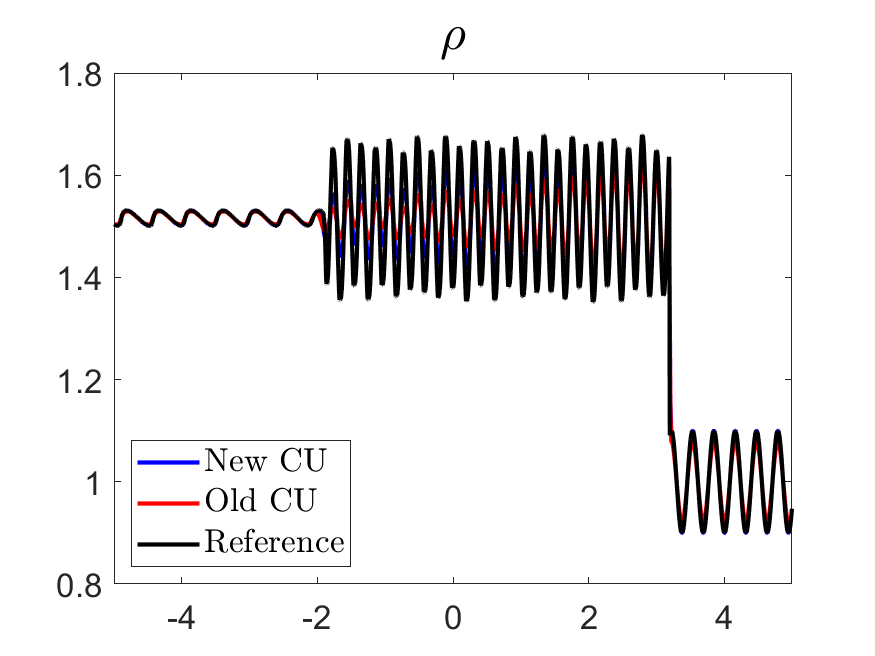}\hspace*{1.0cm}
\includegraphics[trim=0.9cm 0.4cm 1.3cm 0.4cm, clip, width=6.4cm]{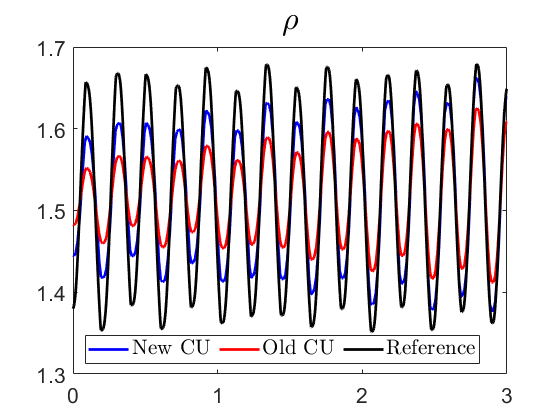}}
\caption{\sf Example 2: Density ($\rho$) computed by the New and Old CU schemes and zoom at $x\in[0,3]$.\label{fig2}}
\end{figure}

In this example, we consider the shock-entropy problem from \cite{Shu88}. The initial conditions,
\begin{equation*}
(\rho,u,p)(x,0)=\begin{cases}
(1.51695,0.523346,1.805),&x<-4.5,\\
(1+0.1\sin(20x),0,1),&x>-4.5,
\end{cases}
\end{equation*}
correspond to a forward-facing shock wave of Mach 1.1 interacting with high-frequency density perturbations, that is, as the shock wave
moves, the perturbations spread ahead. We set free boundary condition at the both ends of the computational domain $[-5,5]$. We apply both
New and Old CU schemes and compute the solutions until the final time $t=5$ on a uniform mesh with $\dx=1/80$. The numerical results are
shown in Figure \ref{fig2} along with the reference solution computed by the Old CU scheme on a much finer mesh with $\dx=1/400$. One can
observe that the New CU scheme produces substantially more accurate results compared to those obtained by the Old CU scheme. This can also
be clearly seen on Figure \ref{fig2} (right), where we zoom at the area where the solution has smooth oscillatory structures.

\subsubsection*{Example 3---Stationary Contact Wave, Traveling Shock and Rarefaction Wave}
\begin{figure}[ht!]
\centerline{
\includegraphics[trim=1.4cm 0.4cm 1.2cm 0.4cm, clip, width=6.4cm]{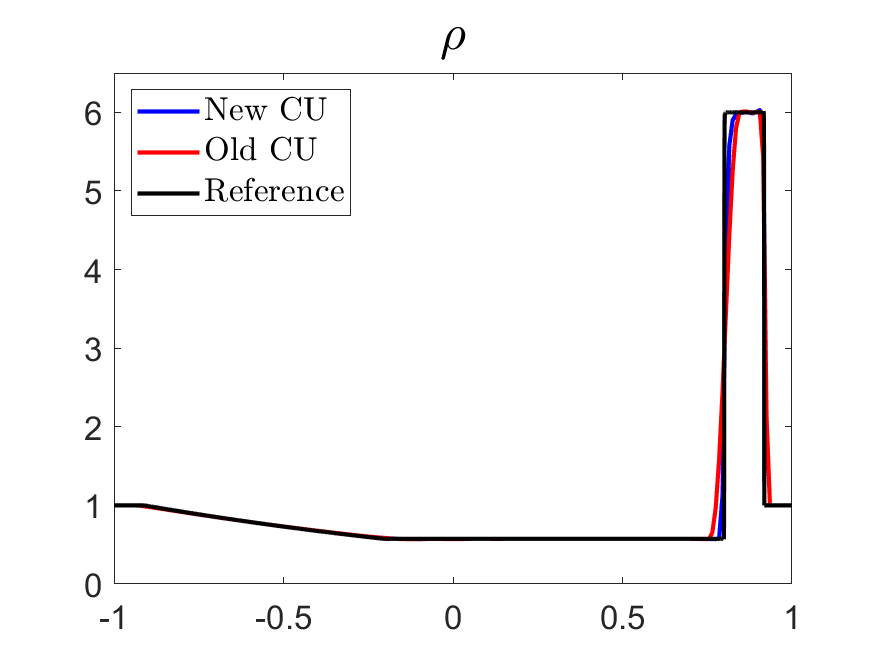}\hspace*{1.0cm}
\includegraphics[trim=1.4cm 0.4cm 1.2cm 0.4cm, clip, width=6.4cm]{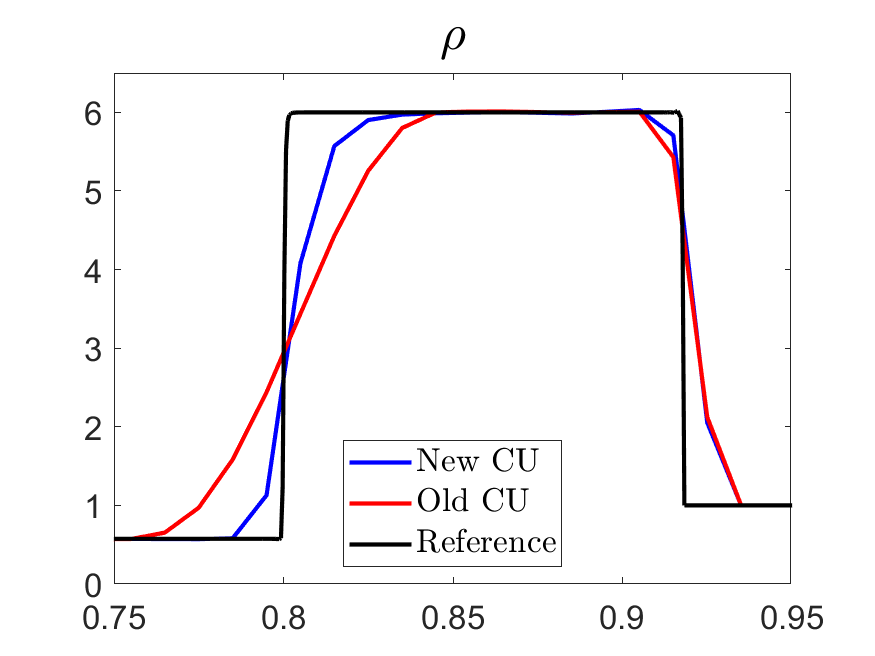}}
\caption{\sf Example 3: Density ($\rho$) computed by the New and Old CU schemes and zoom at $x\in[0.75,0.95]$.\label{fig3}}
\end{figure}

In the third example, we consider the strong shocks interaction problem proposed in \cite{Woodward88}. The initial conditions,
\begin{equation*}
(\rho,u,p)(x,0)=\begin{cases}
(1,-19.59745,1000)&\mbox{if}~x<0.8,\\
(1,-19.59745,0.01)&\mbox{otherwise},\\
\end{cases}
\end{equation*}
are prescribed in the computational domain $[-1,1]$, in which free boundary conditions are implemented. We compute the numerical solutions
until the final time $t=0.03$ by both the New and Old CU schemes on a uniform mesh with $\dx=1/100$, as well as the reference solution,
which is obtained by the Old CU scheme on a much finer mesh with $\dx=1/2000$. The numerical results, plotted in Figure \ref{fig3}, show
that both schemes produce non-oscillatory numerical solutions, but the resolution of the contact wave achieved by the New CU scheme is
higher; see also Figure \ref{fig3} (right), where we zoom at the neighborhood of the contact wave.

\subsubsection*{Example 4---Blast Wave Problem}
\begin{figure}[ht!]
\centerline{
\includegraphics[trim=1.4cm 0.4cm 0.8cm 0.4cm, clip, width=6.5cm]{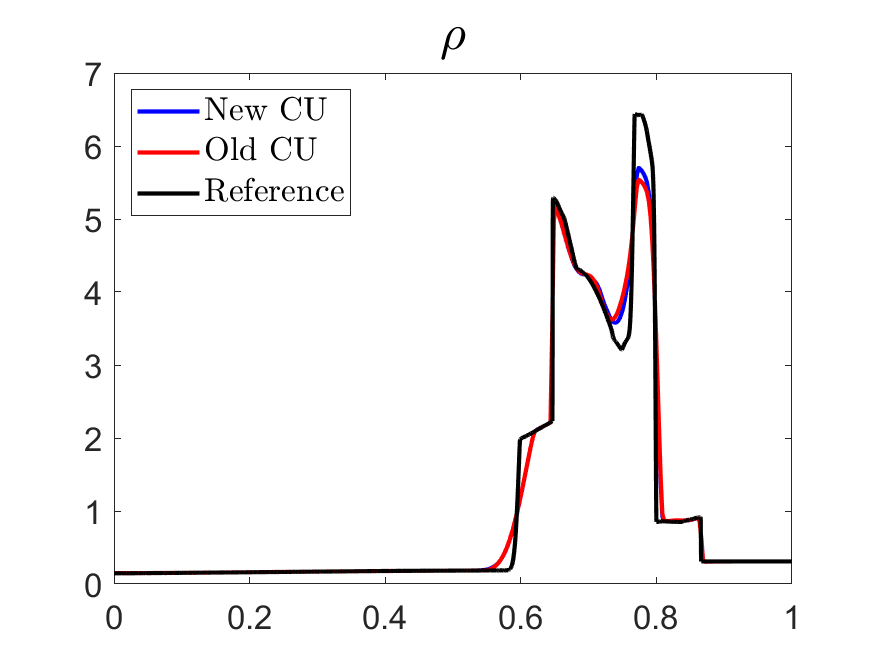}\hspace*{1.0cm}
\includegraphics[trim=1.4cm 0.4cm 0.8cm 0.4cm, clip, width=6.5cm]{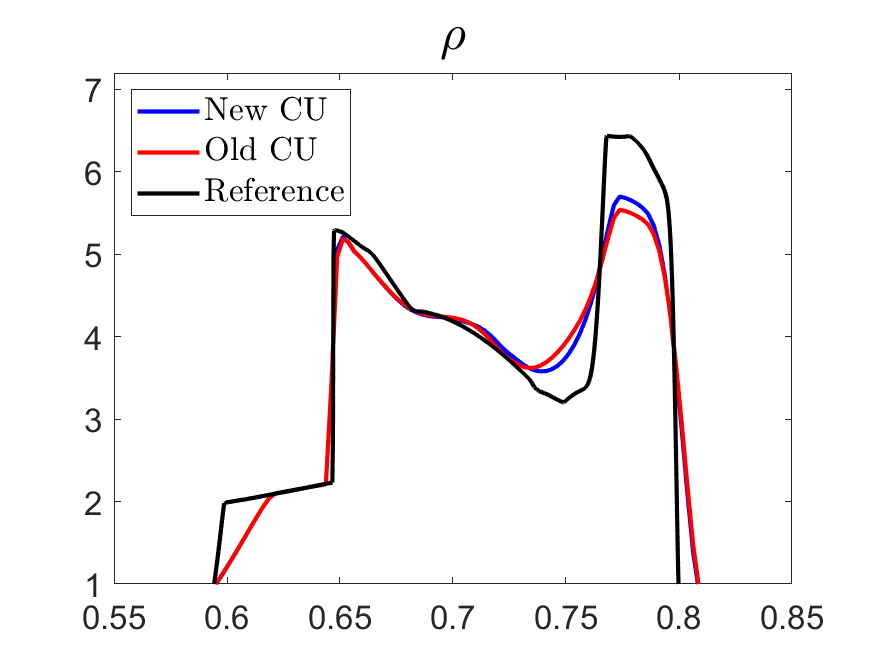}}
\caption{\sf Example 4: Density ($\rho$) computed by the New and Old CU schemes and zoom at $x\in[0.55,0.85]$.\label{fig41}}
\end{figure}

In the final 1-D example, we solve the strong shocks interaction problem from \cite{Woodward88}, which is considered on the interval $[0,1]$
with the solid wall boundary conditions at both ends and subject to the following initial conditions:
\begin{equation*}
(\rho, u,p)(x,0)=\begin{cases}
(1,0,1000),&x<0.1,\\
(1,0,0.01),&0.1\le x\le 0.9,\\
(1,0,100),&x>0.9.
\end{cases}
\end{equation*}
We compute the numerical solutions until the final time $t=0.038$ by both the New and Old CU schemes on a uniform mesh with $\dx=1/400$ and,
as before, implement the Old CU scheme on a fine grid with $\dx=1/4000$ to obtain the corresponding reference solution. The obtained
results, presented in Figure \ref{fig41}, demonstrate that the New CU scheme achieves slightly higher resolution of the second density
spike.

\subsection{Two-Dimensional Examples}
\subsubsection*{Example 5---2-D Riemann Problem}
\begin{figure}[ht!]
\centerline{
\includegraphics[trim=1.3cm 0.4cm 1.7cm 0.4cm, clip, width=6.7cm]{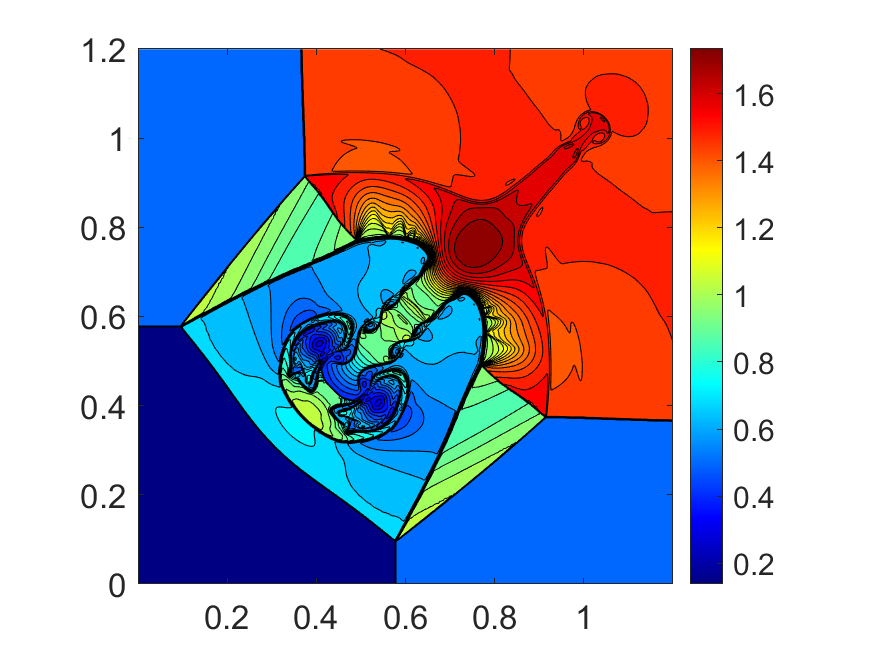}\hspace*{1.0cm}
\includegraphics[trim=1.3cm 0.4cm 1.7cm 0.4cm, clip, width=6.7cm]{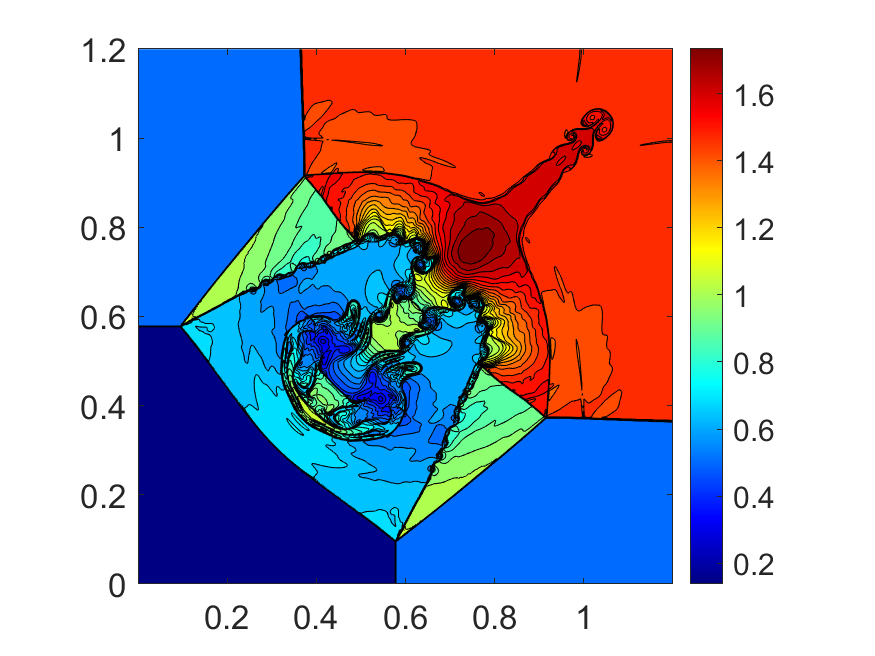}}
\caption{\sf Example 5: Density ($\rho$) computed by the Old (left) and New (right) CU schemes.\label{fig4}}
\end{figure}

In the first 2-D example, we consider Configuration 3 of the 2-D Riemann problems taken from \cite{Kurganov02}; see also
\cite{Schulz93,Schulz93a,Zheng01}. The initial conditions,
\begin{equation*}
(\rho(x,y,0),u(x,y,0),v(x,y,0),p(x,y,0))=\begin{cases}
(1.5,0,0,1.5),&x>1,\,y>1,\\
(0.5323,1.206,0,0.3),&x<1,\,y>1,\\
(0.138, 1.206,1.206,0.029),&x<1,\,y<1,\\
(0.5323,0,1.206,0.3),&x>1,\,y<1,
\end{cases}
\end{equation*}
are prescribed in the computational domain $[0,1.2]\times[0,1.2]$, in which the free boundary conditions are implemented on all of the four
sides of the domain. We compute the numerical solutions until the final time $t=1$ by both the New and Old CU schemes on a uniform mesh with
$\dx=\dy=3/2500$, and plot the obtained results in Figure \ref{fig4}. As one can clearly see, the New CU scheme outperforms the Old CU
scheme in capturing a sideband instability of the jet in the zones of strong along-jet velocity shear and the instability along the jet’s
neck.

\subsubsection*{Example 6---Explosion Problem}
In this example, we consider the explosion problem studied in \cite{Garg21,Kurganov07,Kurganov21a,Liska03}. This is a circularly symmetric
problem with an initial circular region of higher density and pressure with the following initial conditions,
\begin{equation*}
(\rho(x,y,0),u(x,y,0),v(x,y,0),p(x,y,0))=\begin{cases}
(1,0,0,1),&x^2+y^2<0.16,\\
(0.125,0,0,0.1),&\mbox{otherwise},
\end{cases}
\end{equation*}
prescribed in the computational domain $[0,1.5]\times[0,1.5]$. Solid wall boundary conditions are imposed at $x=0$ and $y=0$, while free
boundary conditions are set at $x=1.5$ and $y=1.5$. It is well-known that the solution of this initial-boundary value problem develops
circular shock and contact waves. While the shock wave is stable and a good numerical scheme should contain a sufficient amount of numerical
dissipation to capture the shock in a stable, non-oscillatory manner, the contact wave is unstable and can only be stabilized numerically by
the numerical diffusion present in the scheme. Therefore, this is a good benchmark to measure the amount of numerical dissipation present in
different schemes as one ideally wants to have as little numerical dissipation as possible, but sufficient to stabilize the shock wave.

\begin{figure}[ht!]
\centerline{
\includegraphics[trim=1.3cm 0.4cm 1.4cm 0.4cm, clip, width=6.8cm]{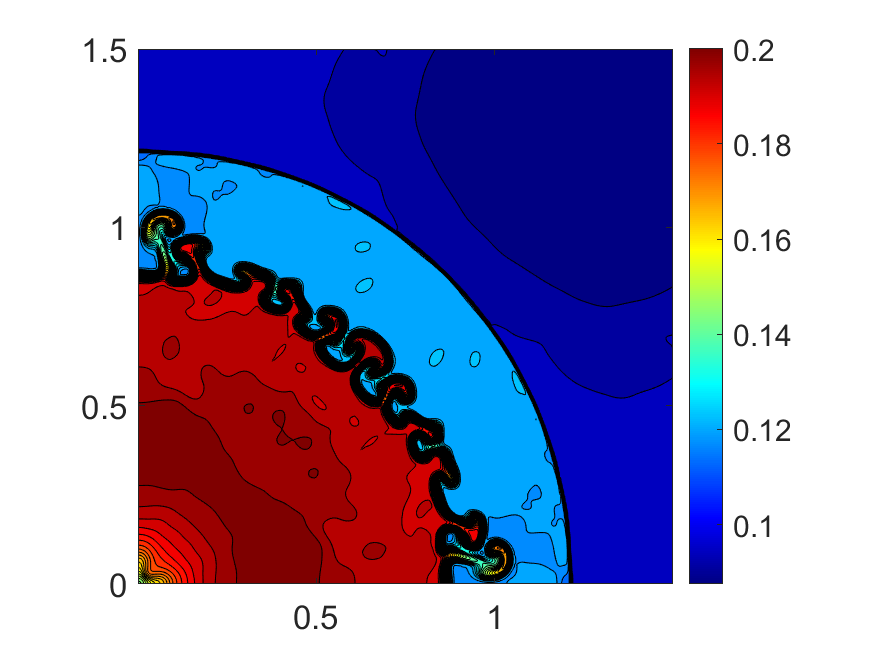}\hspace*{1.0cm}
\includegraphics[trim=1.3cm 0.4cm 1.4cm 0.4cm, clip, width=6.8cm]{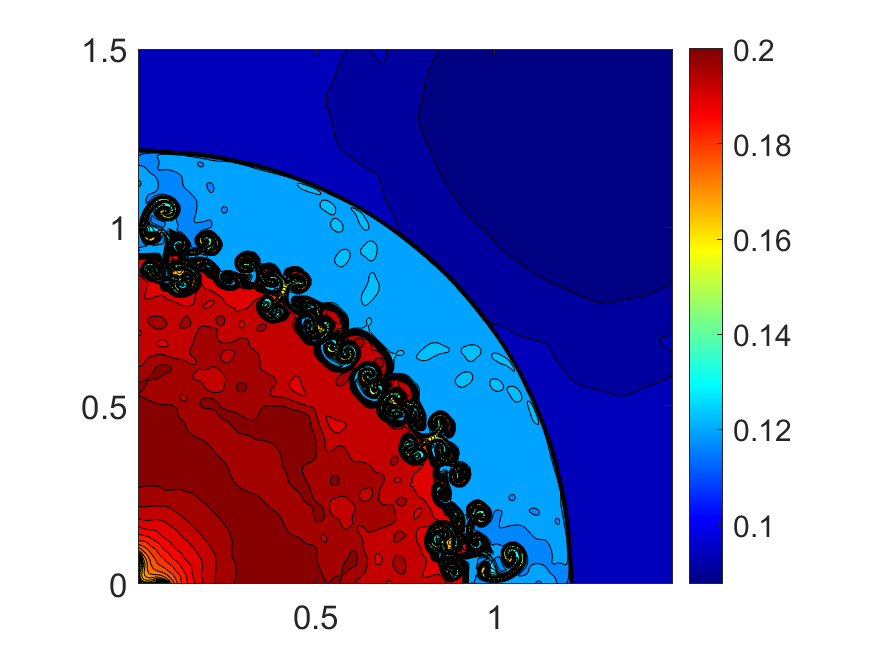}}
\caption{\sf Example 6: Density ($\rho$) computed by the Old (left) and New (right) CU schemes.\label{fig5}}
\end{figure}

We apply both New and Old CU schemes and compute the numerical solutions on a uniform mesh with $\dx=\dy=3/800$ until the final time
$t=3.2$. The obtained results are presented in Figure \ref{fig5}. One can observe that compared with the results obtained by the Old CU
scheme, the contact curve captured by the New CU scheme is substantially ``curlier'' and the mixing layer is slightly wider (indicating a
more severe instability), while the shock is still stable.

\subsubsection*{Example 7---Implosion Problem}
\begin{figure}[ht!]
\centerline{
\includegraphics[trim=1.3cm 0.4cm 1.7cm 0.4cm, clip, width=6.7cm]{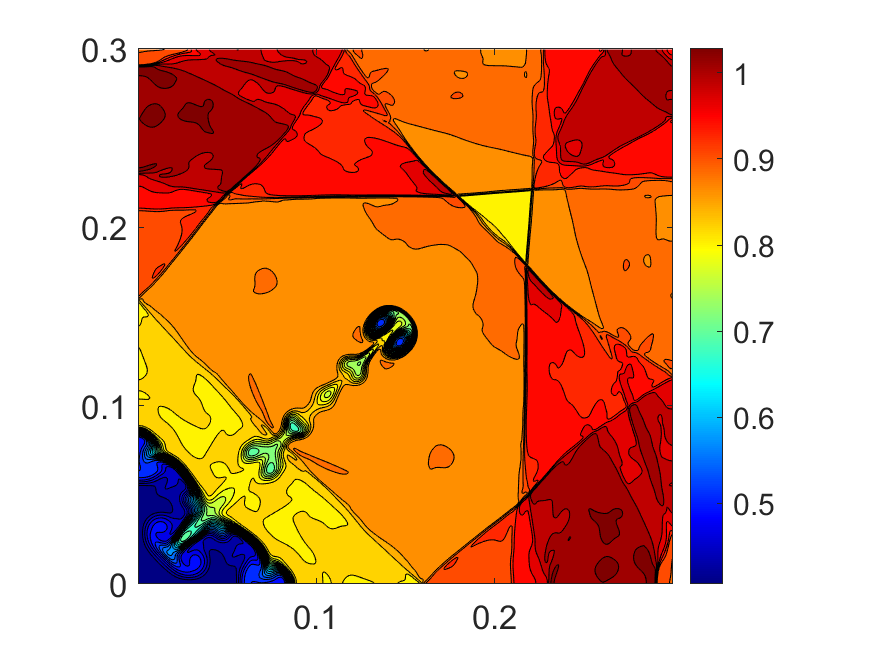}\hspace*{1.0cm}
\includegraphics[trim=1.3cm 0.4cm 1.7cm 0.4cm, clip, width=6.7cm]{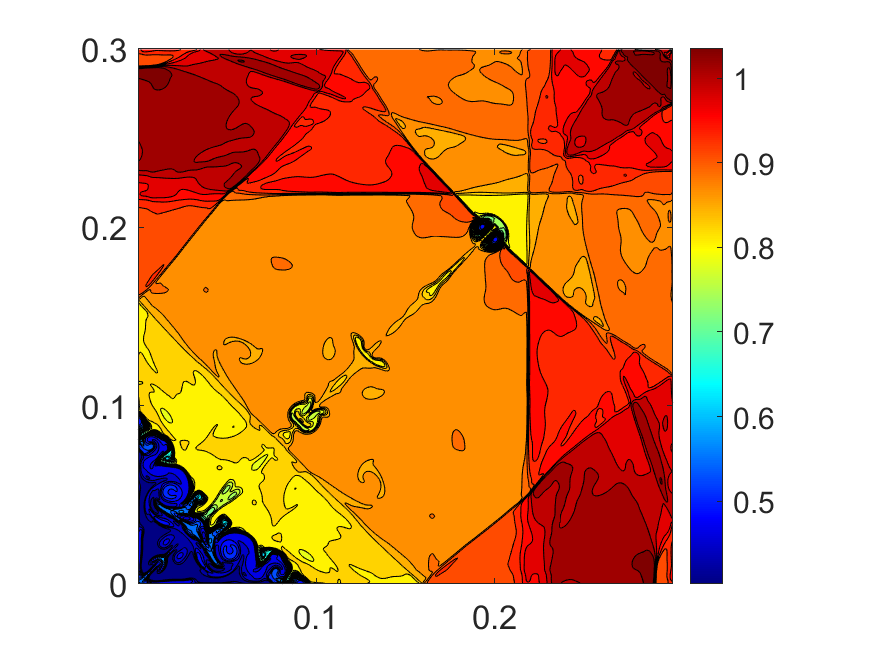}}
\vskip5pt
\centerline{
\includegraphics[trim=1.3cm 0.4cm 1.7cm 0.4cm, clip, width=6.7cm]{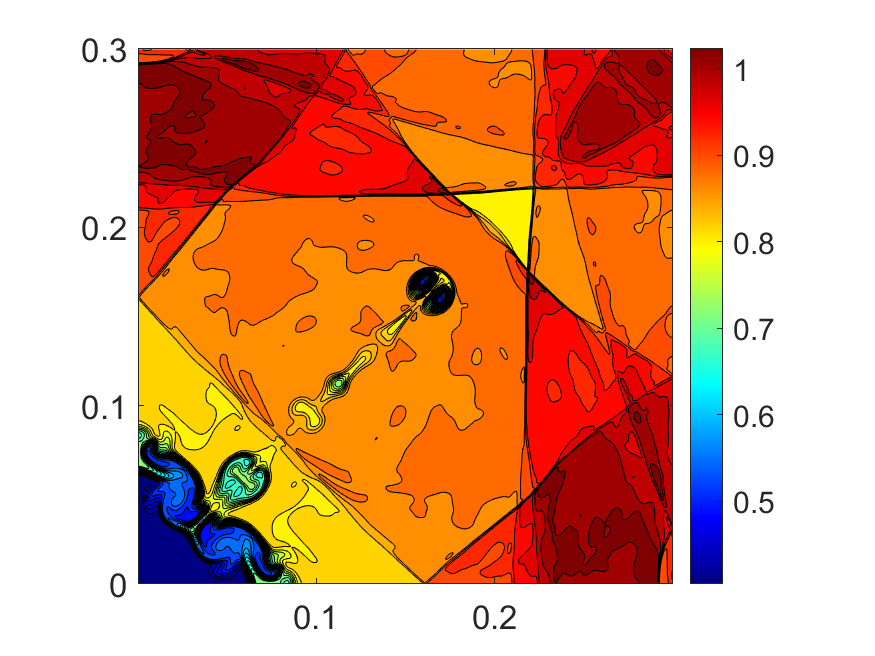}\hspace*{1.0cm}
\includegraphics[trim=1.3cm 0.4cm 1.7cm 0.4cm, clip, width=6.7cm]{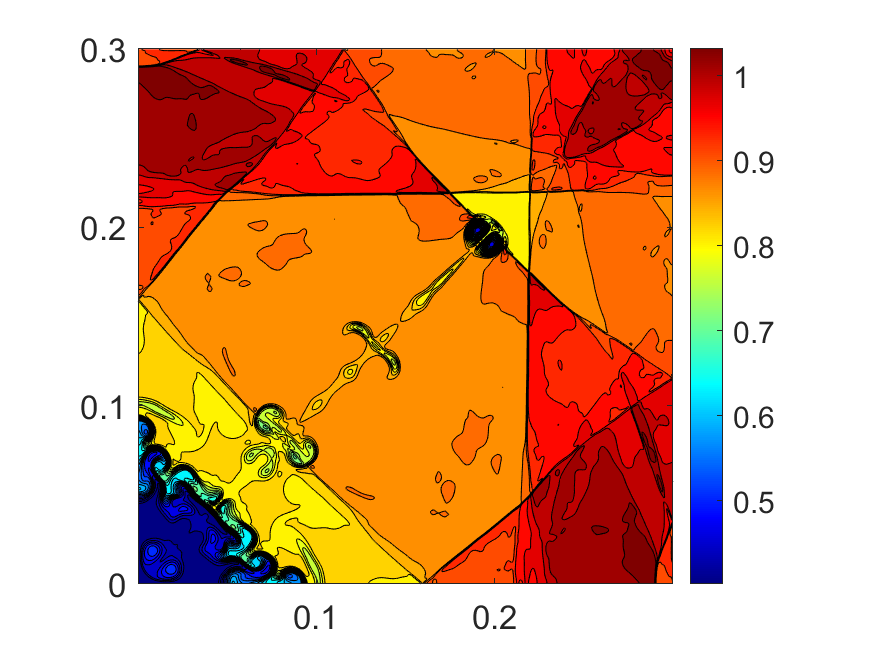}}
\caption{\sf{Example 7: Density ($\rho$) computed by the Old CU scheme using $\dx=\dy=1/2000$ (top left), New CU schemes using
$\dx=\dy=1/2000$ (top right), Old CU scheme using $\dx=\dy=1/2500$ (bottom left), and  Old CU scheme using $\dx=\dy=3/10000$ (bottom
right).}\label{fig6}}
\end{figure}

In this example, we consider the implosion problem taken from \cite{Garg21,Kurganov07,Kurganov21a,Liska03}. The initial conditions,
\begin{equation*}
(\rho(x,y,0),u(x,y,0),v(x,y,0),p(x,y,0))=\begin{cases}
(0.125,0,0,0.14),&|x|+|y|<0.15,\\
(1,0,0,1),&\mbox{otherwise},
\end{cases}
\end{equation*}
are prescribed in the computational domain $[0,0.3]\times[0,0.3]$ with solid boundary conditions imposed at all of the four sides. We
compute the numerical solutions until the final time $t=2.5$ by both the New and Old CU schemes on a uniform mesh with $\dx=\dy=1/2000$. The
obtained results are depicted in Figure \ref{fig6} (top row). As one can clearly see, a jet generated by the New CU scheme propagates much
further in the direction of $y=x$ than the jet produced by the Old CU scheme: This is attributed to a much smaller amount of numerical
dissipation present in the New CU scheme.

In this example, we also check the efficiency of the proposed LCD-based CU scheme by performing a more thorough comparison between the New
and Old CU schemes: We now take into account additional computational cost of the New CU scheme. To this end, we measure the CPU time
consumed during the above calculations by the New CU scheme and refine the mesh used by the Old CU scheme to the level that exactly the same
CPU time is consumed to compute both of the numerical solutions. The corresponding grids are $\dx=\dy=1/2000$ for the New CU scheme, and
$\dx=\dy=1/2500$ for the Old CU scheme (this solution is plotted in Figure \ref{fig6} (bottom left). The obtained numerical results indicate
that the New CU scheme still achieves a much higher resolution than the Old CU scheme.

It is instructive to check whether the solution computed by the Old CU scheme on even finer mesh will be comparable with the New CU
solution.  To this end, we refine the mesh to the level $\dx=\dy=3/10000$ and plot the obtained Old CU scheme in Figure \ref{fig6} (bottom
right). As one can see, this solution is similar to the New CU solution shown in Figure \ref{fig6} (top right). This shows that the New CU
scheme is capable of resolving the same details on a substantially coarser grid due to the smaller amount of numerical diffusion present in
the scheme.

\subsubsection*{Example 8---Kelvin-Helmholtz (KH) Instability}
\begin{figure}[ht!]
\centerline{
\includegraphics[trim=0.6cm 0.4cm 1.1cm 0.3cm, clip, width=6.7cm]{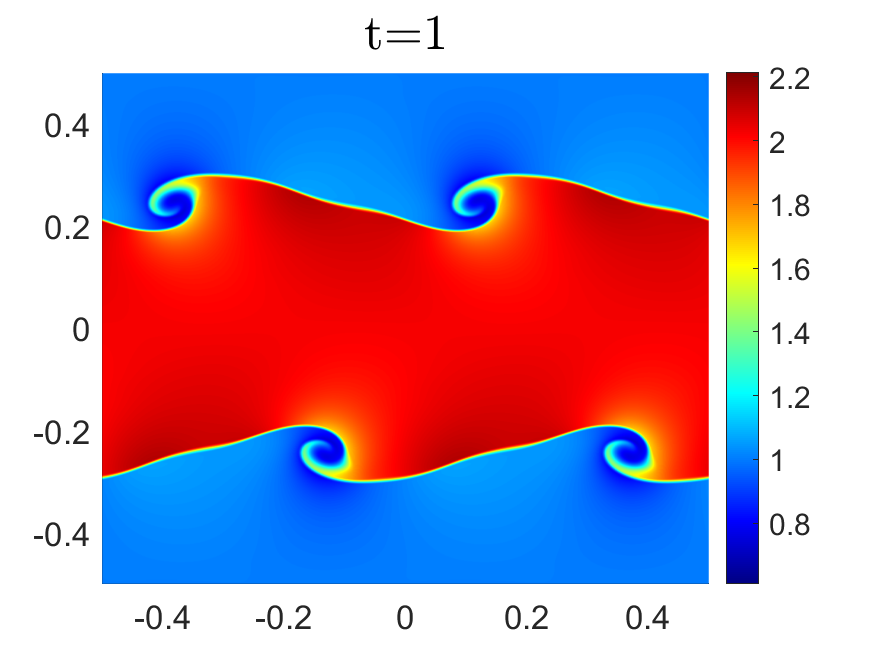}\hspace*{1.0cm}
\includegraphics[trim=0.6cm 0.4cm 1.1cm 0.3cm, clip, width=6.7cm]{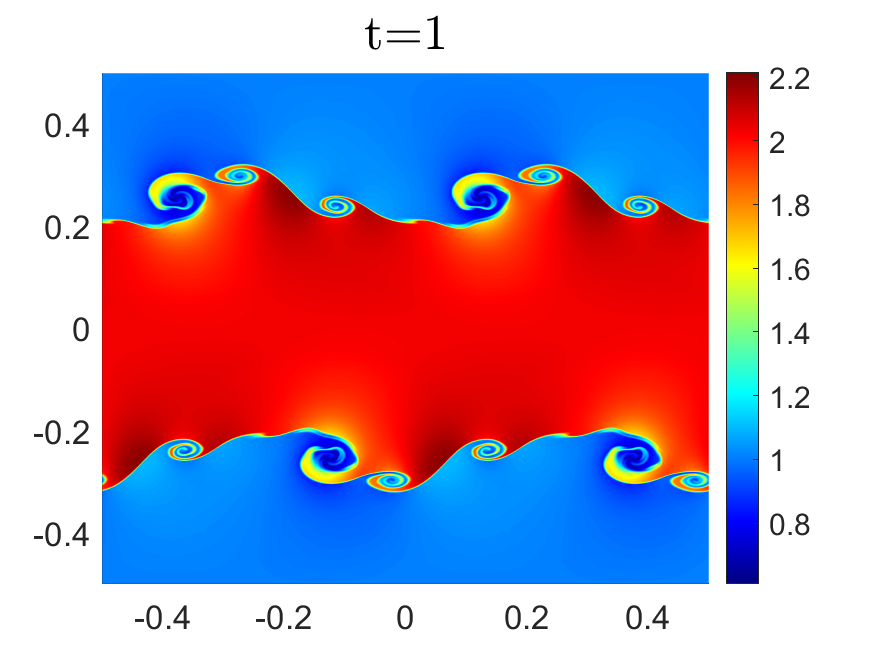}}
\vskip5pt
\centerline{
\includegraphics[trim=0.6cm 0.4cm 1.1cm 0.3cm, clip, width=6.7cm]{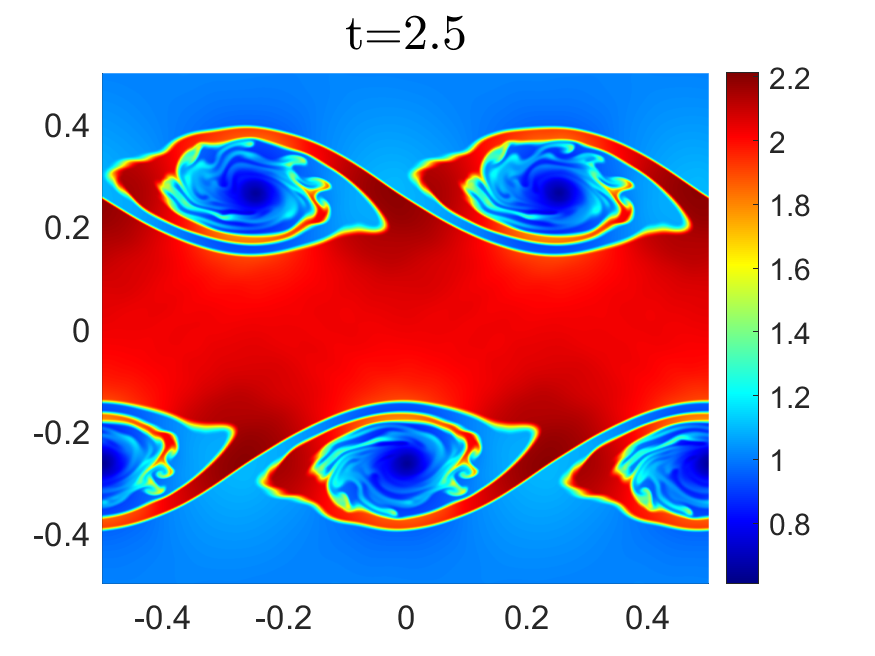}\hspace*{1.0cm}
\includegraphics[trim=0.6cm 0.4cm 1.1cm 0.3cm, clip, width=6.7cm]{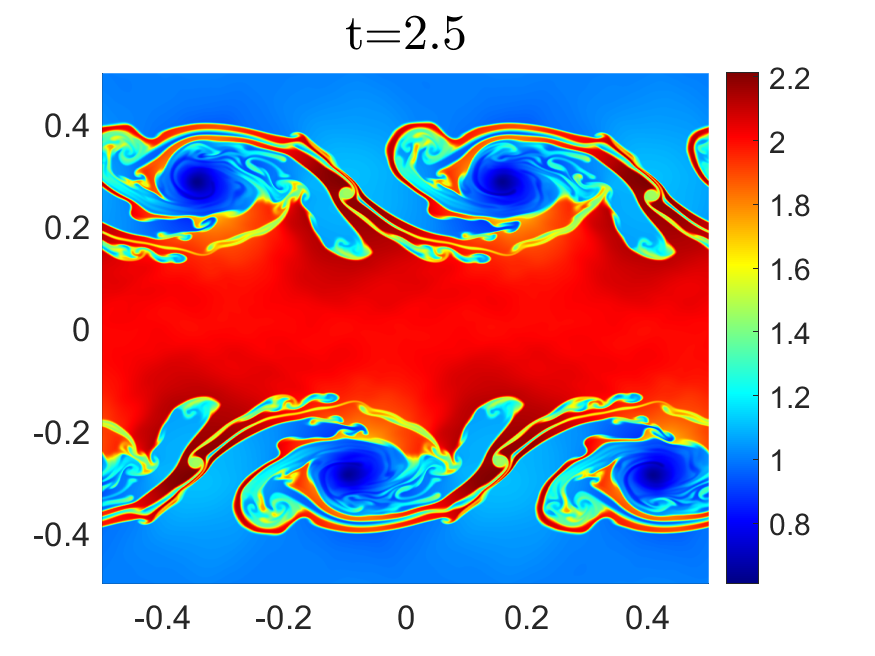}}
\vskip5pt
\centerline{
\includegraphics[trim=0.6cm 0.4cm 1.1cm 0.3cm, clip, width=6.7cm]{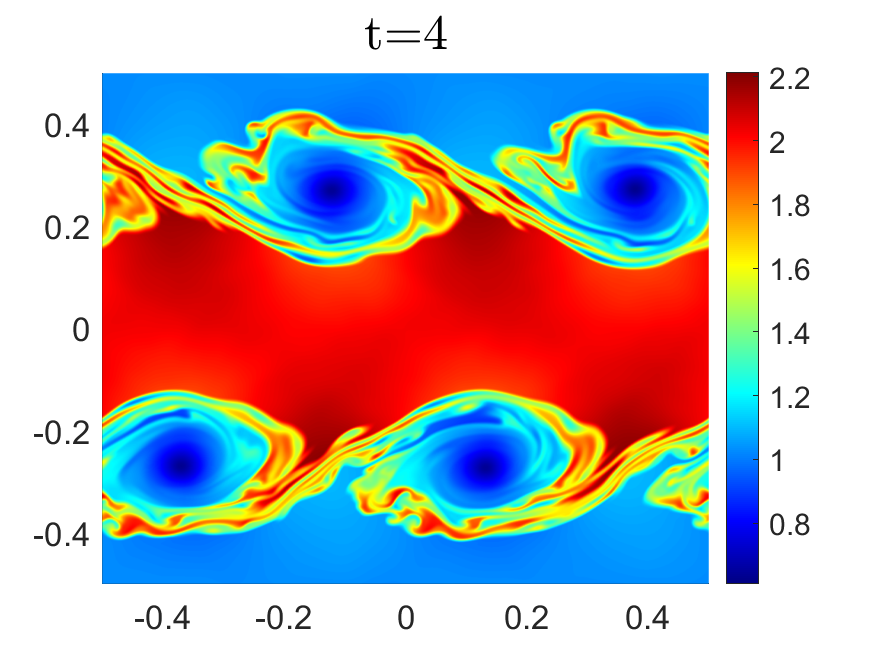}\hspace*{1.0cm}
\includegraphics[trim=0.6cm 0.4cm 1.1cm 0.3cm, clip, width=6.7cm]{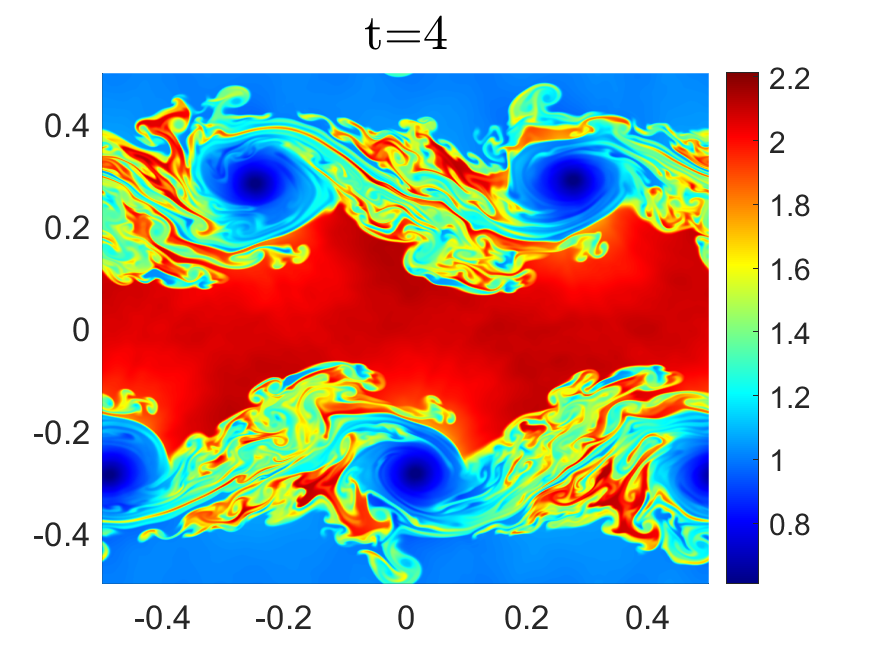}}
\caption{\sf Example 8: Time snapshots of density ($\rho$) computed by the Old (left column) and New (right column) CU schemes at $t=1$ (top
row), $t=2.5$ (middle row), and $t=4$ (bottom row).\label{fig7}}
\end{figure}

In this example, we study the KH instability, which develops in the test problem taken from
\cite{Feireisl21,Fjordholm16,Garg21,Kurganov21a,Panuelos20}. We take the following initial data:
\begin{equation*}
\begin{aligned}
&(\rho(x,y,0),u(x,y,0))=\begin{cases}
(1,-0.5+0.5e^{(y+0.25)/L}),&y<-0.25,\\
(2,0.5-0.5e^{(-y-0.25)/L}),&-0.25<y<0,\\
(2,0.5-0.5e^{(y-0.25)/L}),&0<y<0.25,\\
(1,-0.5+0.5e^{(0.25-y)/L}),&y>0.25,
\end{cases}\\
&v(x,y,0)=0.01\sin(4\pi x),\qquad p(x,y,0)\equiv1.5,
\end{aligned}
\end{equation*}
where $L$ is a smoothing parameter (we take $L=0.00625$), which corresponds to a thin shear interface with a perturbed vertical velocity
field $v$ in the conducted simulations. Periodic boundary conditions are imposed on all of the four sides of the computational domain
$[-0.5,0.5]\times[-0.5,0.5]$. We compute the numerical solutions until the final time $t=4$ by both the New and Old CU schemes on a uniform
mesh with $\dx=\dy=1/1024$. The numerical results at times $t=1$, 2.5, and 4 are presented in Figure \ref{fig7}. As one can see, at the
early time $t=1$, the vortex streets formed by the adaptive scheme are more pronounced. These structures continue growing exponentially in
time showing much more complicated turbulent mixing captured by the New CU scheme, especially at the final time $t=4$. Clearly, the New CU
scheme outperforms the Old one in achieving a higher resolution of the KH instabilities.

\begin{figure}[ht!]
\centerline{
\includegraphics[trim=0.6cm 0.4cm 1.1cm 0.1cm, clip, width=6.8cm]{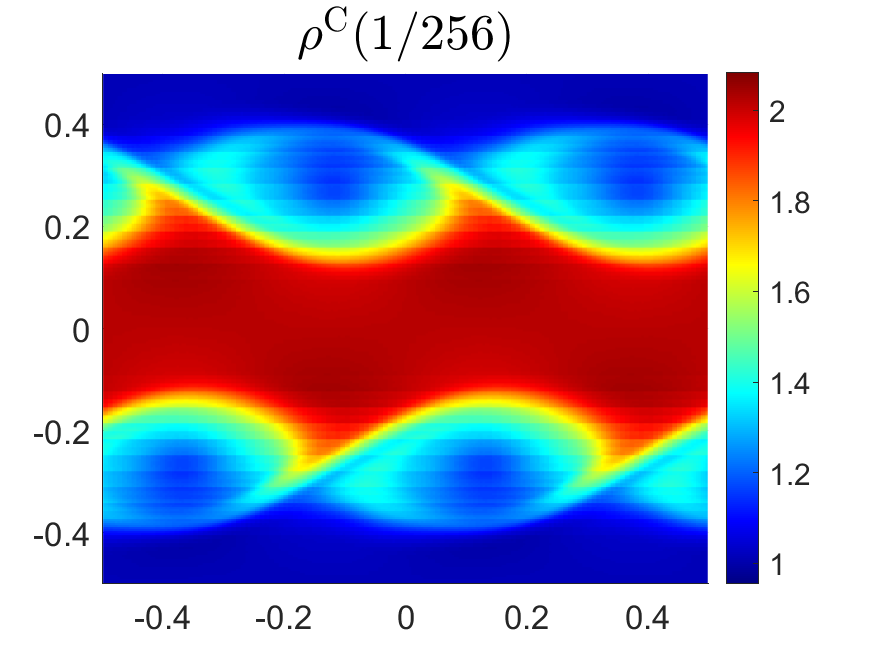}\hspace*{1.0cm}
\includegraphics[trim=0.6cm 0.4cm 1.1cm 0.1cm, clip, width=6.8cm]{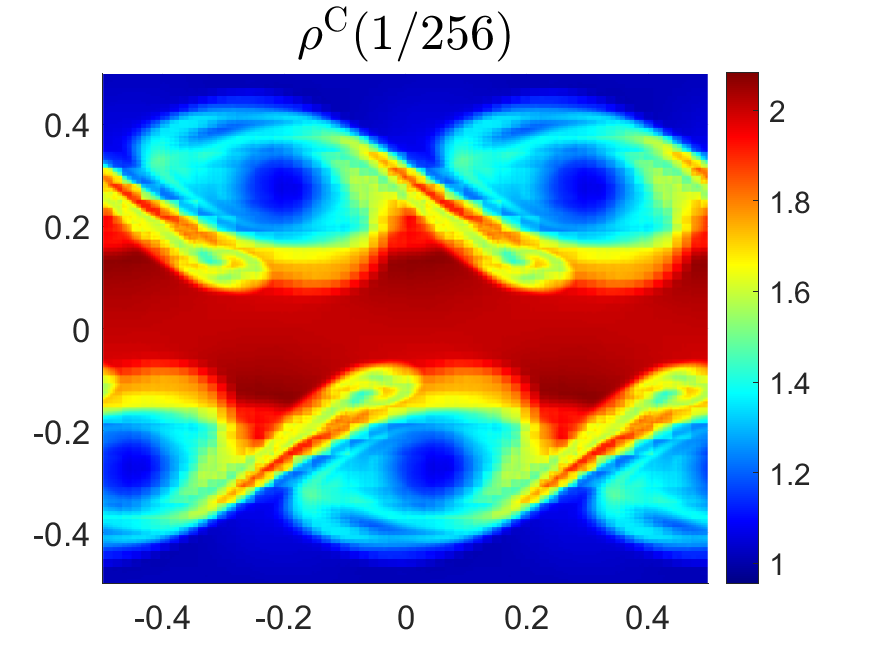}}
\vskip5pt
\centerline{
\includegraphics[trim=0.6cm 0.4cm 1.1cm 0.1cm, clip, width=6.8cm]{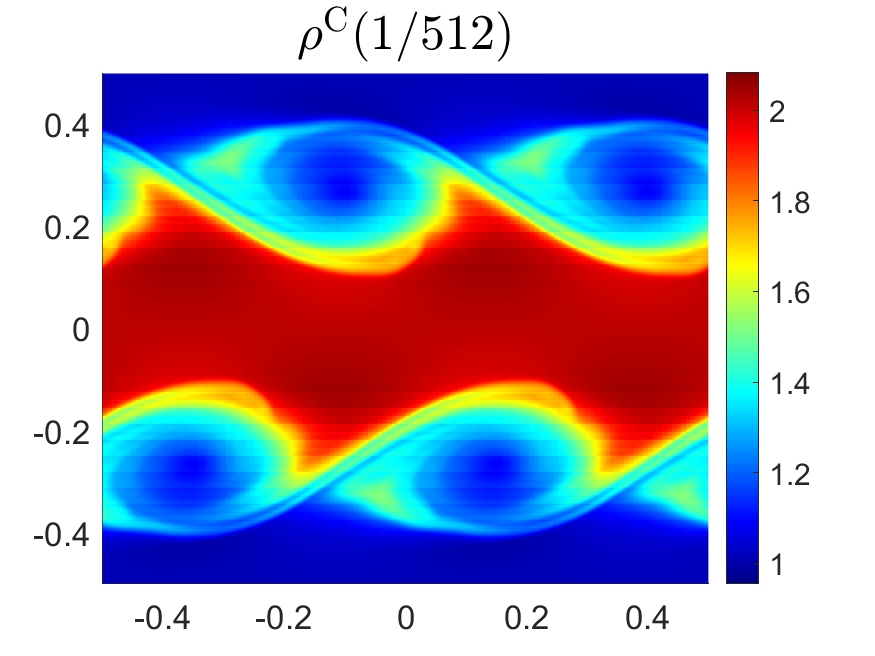}\hspace*{1.0cm}
\includegraphics[trim=0.6cm 0.4cm 1.1cm 0.1cm, clip, width=6.8cm]{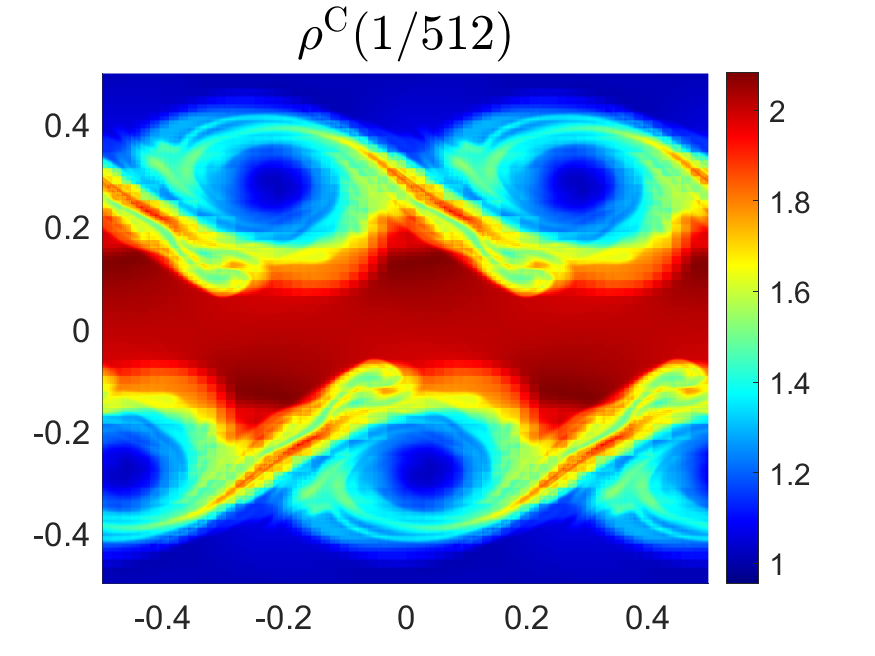}}
\vskip5pt
\centerline{
\includegraphics[trim=0.6cm 0.4cm 1.1cm 0.1cm, clip, width=6.8cm]{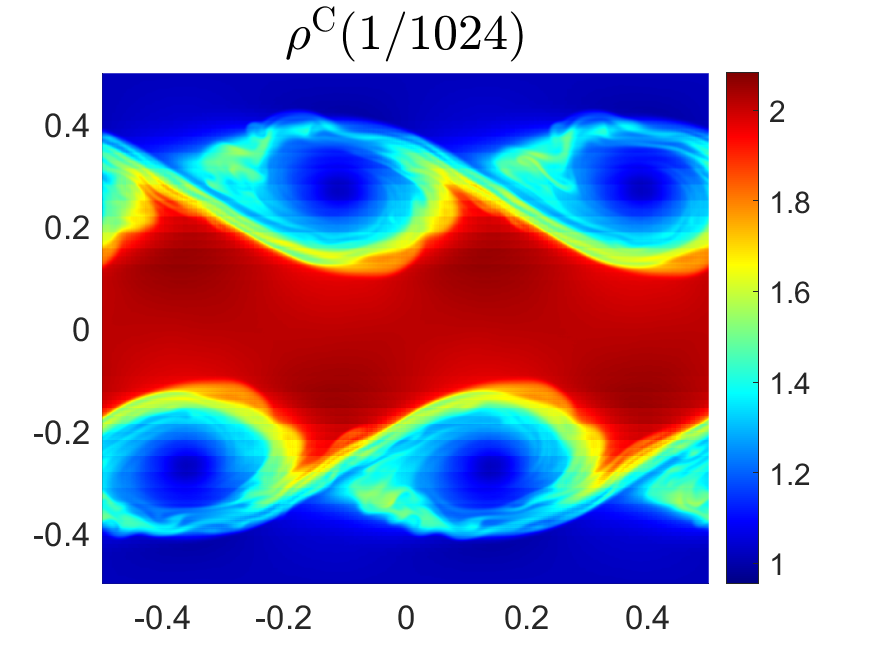}\hspace*{1.0cm}
\includegraphics[trim=0.6cm 0.4cm 1.1cm 0.1cm, clip, width=6.8cm]{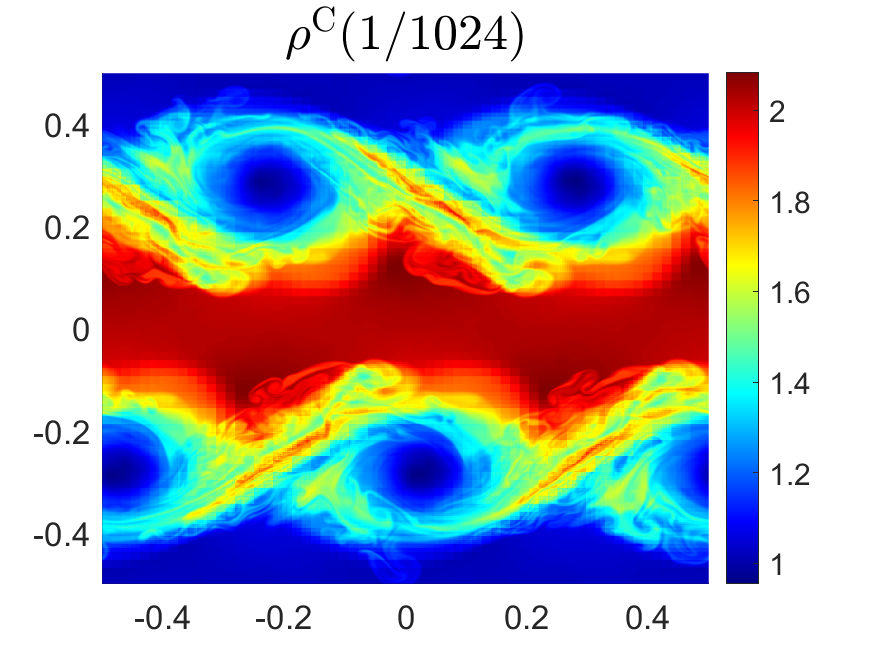}}
\caption{\sf Example 8: Ces\`aro averages of the density $\rho^{\rm C}$ computed by the Old (left column) and New (right column) CU schemes
at $t=4$.\label{fig71}}
\end{figure}

For the KH instability problem, it is well-known that numerical solutions do not converge when the mesh is refined, since the limiting
solution is not a weak solution but a so-called dissipative weak solution; see \cite{Lukacova_book} for more details. The latter can be
approximated by the Ces\`aro averages computed from the numerical solutions on meshes with different resolutions. As proved in
\cite{Lukacova_book}, the Ces\`aro averages converge strongly (in the $L^q$-norm with $1\le q<\infty$) to a dissipative weak solution. In
order to approximate the limiting solution, we consider the Ces\`aro averages of the densities computed at the final time $t=4$ by the New
and Old CU schemes. To this end, we introduce the sequence of meshes with the cells of size ${1}/{2^n}$, $n=5,\ldots,10$, and denote by
$\rho(1/2^n)$ the density computed on the corresponding mesh. We then project the obtained coarser mesh solutions with $n=5,\ldots,m-1$ onto
the finer mesh with $n=m$ (the projection is carried out using the minmod reconstruction \cite{Lie03,Nessyahu90,Sweby84} applied to the
density field) and denote the obtained results still by $\rho(1/2^n)$, $n=5,\ldots,m$. After this, we compute the Ces\`aro averages by
$$
\rho^{\rm C}(1/2^m)=\frac{\rho(1/2^5)+\cdots+\rho(1/2^n)}{m-4},\quad m=8,9,10.
$$
The obtained results, presented in Figure \ref{fig71} for both the New and Old CU schemes, indicate non-uniqueness of the limiting
dissipative weak solution that apparently depends on the choice of the numerical diffusion. This leads to an interesting question on
suitable selection criteria that can be studied in the future.

\subsubsection*{Example 9---Rayleigh-Taylor (RT) Instability}
In the last example taken from \cite{Garg21,Kurganov21a,Shi03,Wang20}, we investigate the RT instability, which is a physical phenomenon
occurring when a layer of heavier fluid is placed on top of a layer of lighter fluid. The model is governed by the 2-D Euler equations
\eref{1.1}, \eref{6.3}--\eref{6.4} with added gravitational source terms. In the studied setup, the gravitational force acts in the
positive $y$-direction and the modified system reads as
$$
\begin{aligned}
&\rho_t+(\rho u)_x+(\rho v)_y=0,\\
&(\rho u)_t+(\rho u^2 +p)_x+(\rho uv)_y=0,\\
&(\rho v)_t+(\rho uv)_x+(\rho v^2+p)_y=\rho,\\
&E_t+[u(E+p)]_x+[v(E+p)]_y=\rho v.
\end{aligned}
$$
This system is considered subject to the following initial conditions:
\begin{equation*}
(\rho(x,y,0),u(x,y,0),v(x,y,0),p(x,y,0))=\begin{cases}
(2,0,-0.025\,c\cos(8\pi x),2y+1),&y<0.5,\\
(1,0,-0.025\,c\cos(8\pi x),y+1.5),&\mbox{otherwise},
\end{cases}
\end{equation*}
where $c:=\sqrt{\gamma p/\rho}$ is the speed of sound, and the solid wall boundary conditions at $x=0$ and $x=0.25$, and the following
Dirichlet boundary conditions at the top and bottom boundaries:
$$
(\rho,u,v,p)|_{y=1}=(1,0,0,2.5),\qquad(\rho,u,v,p)|_{y=0}=(2,0,0,1).
$$

\begin{figure}[ht!]
\centerline{\includegraphics[width=17.8cm]{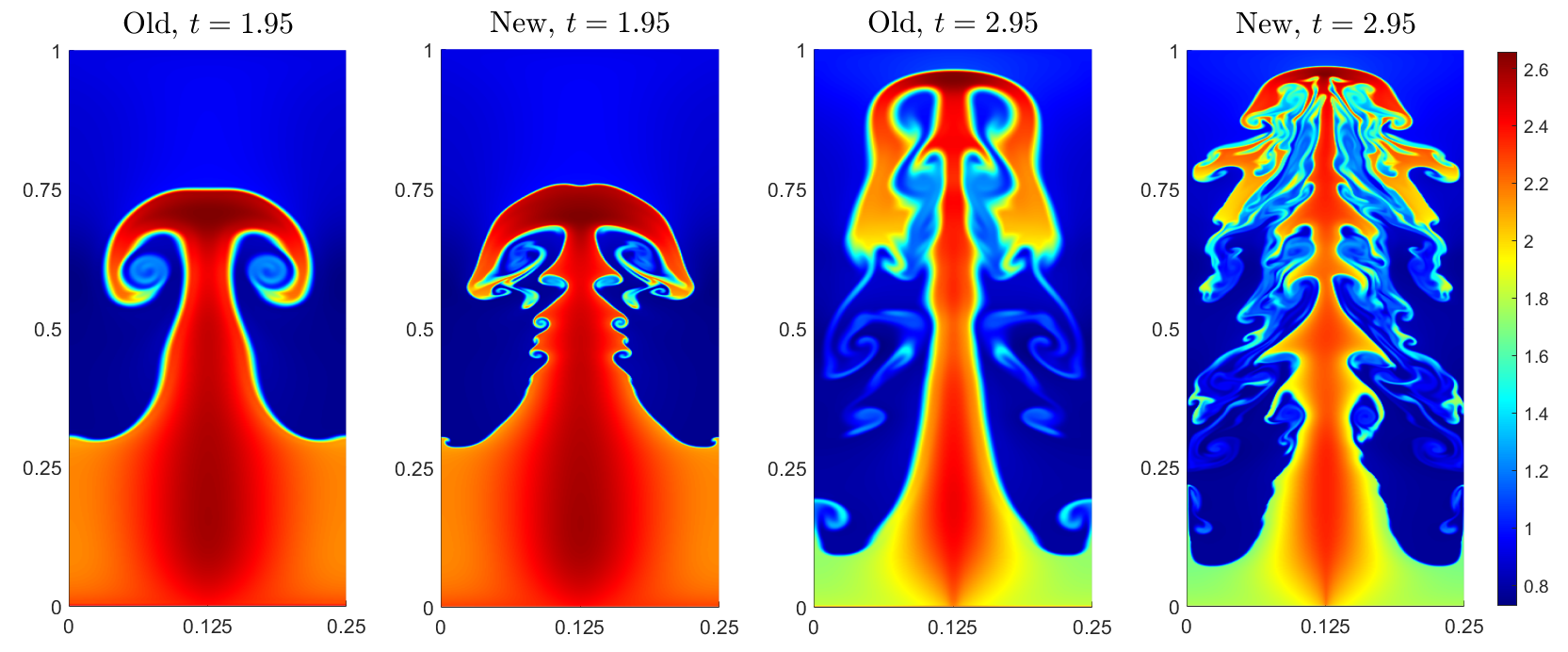}}
\caption{\sf Example 9: Density ($\rho$) computed by the Old and New CU schemes at different times.\label{fig8}}
\end{figure}

We compute the numerical solutions in the computational domain $[0,0.25]\times[0,1]$ discretized using a uniform mesh with $\dx=\dy=1/1024$
until the final time $t=2.95$ by both the New and Old CU schemes. The numerical results at times $t=1.95$ and 2.95, presented in Figure
\ref{fig8}, show a significant difference in performance of the New and Old CU schemes. Indeed, the New CU scheme is capable of capturing the
RT instability with a much higher resolution of the complicated solution structure. Once again, this indicates that the New CU scheme is
substantially less dissipative than its Old counterpart.

\subsection*{Acknowledgment}
The work of A. Chertock was supported in part by NSF grant DMS-1818684. The work of M. Herty was supported in part by the DFG (German
Research Foundation) through 20021702/GRK2326, 333849990/IRTG-2379, HE5386/18-1, 19-2, 22-1, 23-1 and under Germany’s Excellence Strategy
EXC-2023 Internet of Production 390621612. The work of A. Kurganov was supported in part by NSFC grants 12171226 and 12111530004, and by the
fund of the Guangdong Provincial Key Laboratory of Computational Science and Material Design (No. 2019B030301001). The work of M.
Luk\'a\v{c}ov\'a-Medvi{\softd}ov\'a has been funded by the DFG (German Research Foundation) under the SFB/TRR 146 Multiscale Simulation
Methods for Soft Matter Systems. M. Luk\'a\v{c}ov\'a-Medvi{\softd}ov\'a gratefully acknowledges support of the Gutenberg Research College of
University Mainz and the Mainz Institute of Multiscale Modeling.

\appendix
\section{LCD Matrices for the 1-D Euler Equations of Gas Dynamics}\label{appa}
The 1-D Euler equations of gas dynamics read as the system \eref{2.1} with
$$
\mU=\begin{pmatrix}\rho\\\rho u\\E\end{pmatrix}\quad\mbox{and}\quad\mF(\mU)=\begin{pmatrix}\rho u\\\rho u^2+p\\u(E+p)\end{pmatrix},
$$
where $\rho$, $u$, $p$, and $E$ are the density, velocity, pressure and total energy, respectively. The system is completed through the
following equations of state:
\begin{equation*}
E=\rho e+\hf\rho u^2,\quad e=\frac{p}{(\gamma-1)\rho},
\end{equation*}
where $e$ is the specific internal energy, and $\gamma$ is the parameter representing the specific heat ratio.

We first compute the Jacobian
$$
\begin{aligned}
A(\mU)=\frac{\partial\mF}{\partial\mU}(\mU)=\begin{pmatrix}0&1&0\\
\dfrac{\gamma-3}{2}u^2&(3-\gamma)u&\gamma-1\\[1.5ex]
-\dfrac{\gamma u\hat E}{\rho}+(\gamma-1)u^3&H-(\gamma-1)u^2&\gamma u\end{pmatrix},
\end{aligned}
$$
where $H=\frac{E+p}{\rho}$ is the total specific enthalpy. We then introduce the matrices
$$
\begin{aligned}
\widehat A_\jph=\begin{pmatrix}0&1&0\\
\dfrac{\gamma-3}{2}\hat u^2&(3-\gamma)\hat u&\gamma-1\\[1.5ex]
-\dfrac{\gamma\hat u\hat E}{\hat\rho}+(\gamma-1)\hat u^3&\hat H-(\gamma-1)\hat u^2&\gamma\hat u\end{pmatrix},
\end{aligned}
$$
where $\hat{(\cdot)}$ stands for the following averages:
\begin{equation*}
\hat\rho=\frac{\xbar\rho_j+\xbar\rho_{j+1}}{2},~\hat u=\frac{u_j+u_{j+1}}{2},~\hat p=\frac{p_j+p_{j+1}}{2},~
\hat H=\frac{\hat E+\hat p}{\hat\rho},~\hat E=\frac{\hat p}{\gamma-1}+\hf\hat\rho\hat u^2,
\end{equation*}
where
$$
u_j=\frac{(\xbar{\rho u})_j}{\xbar\rho_j}\quad\mbox{and}\quad p_j=(\gamma-1)\left[\xbar E_j-\hf\rho_j u_j^2\right].
$$

Notice that all of the $\hat{(\cdot)}$ quantities have to have a subscript index, that is, $\hat{(\cdot)}=\hat{(\cdot)}_\jph$, but we omit
it for the sake of brevity for all of the quantities except for $\widehat A$. We then compute the matrix $R_\jph$ composed of the right
eigenvectors of $\widehat A_\jph$ and obtain
$$
R_\jph=\begin{pmatrix}1&1&1\\\hat u-\hat c&\hat u&\hat u+\hat c\\\hat H-\hat u\hat c&\dfrac{\hat u^2}{2}&\hat H+\hat u\hat c\end{pmatrix}
\quad\mbox{and}\quad R^{-1}_\jph=\dfrac{1}{\hat\phi}\begin{pmatrix}
\dfrac{\hat u^2}{2}+\dfrac{\hat u\hat\phi}{2\hat c}&-\hat u-\dfrac{\hat\phi}{2\hat c}&1\\[1.2ex]2\hat\phi-{2\hat H}&{2\hat u}&-2\\
\dfrac{\hat u^2}{2}-\dfrac{\hat u\hat\phi}{2\hat c}&-\hat u+\dfrac{\hat\phi}{2\hat c}&1\end{pmatrix},
$$
where $\hat\phi=2\hat H-\hat u^2$ and $\hat c=\sqrt{\gamma\hat p/\hat\rho}$.

\section{LCD-Based Piecewise Linear Reconstruction for the 2-D Euler Equations of Gas Dynamics}\label{appb}
The 2-D Euler equations of gas dynamics read as \eref{1.1} with
\begin{equation}
\mU=\begin{pmatrix}\rho\\\rho u\\\rho v\\E\end{pmatrix},\quad\mF(\mU)=\begin{pmatrix}\rho u\\\rho u^2+p\\\rho uv\\u(E+p)\end{pmatrix},\quad
\mbox{and}\quad\mG(\mU)=\begin{pmatrix}\rho v\\\rho uv\\\rho v^2+p\\v(E+p)\end{pmatrix},
\label{6.3}
\end{equation}
where $\rho$ is the density, $u$ and $v$ are the $x$- and $y$-velocities, $p$ is the pressure, and $E$ is the total energy. The system is
completed through the following equations of state:
\begin{equation}
E=\rho e+\frac{\rho}{2}(u^2+v^2),\quad e=\frac{p}{(\gamma-1)\rho}.
\label{6.4}
\end{equation}
where, as in the 1-D case, $e$ is the specific internal energy, and $\gamma$ is the parameter representing the specific heat ratio.

In what follows, we discuss the reconstruction of the point values $\mU^{\rm E,W}_{j,k}$ and $\bm\Gamma^{\rm E,W}_{j,k}$ (the point values
$\mU^{\rm N,S}_{j,k}$ and $\bm\Gamma^{\rm N,S}_{j,k}$ can be computed in a similar manner and we omit the details for the sake of brevity).
To this end, we first compute the Jacobian
$$
\begin{aligned}
A(\mU)=\frac{\partial\mF}{\partial\mU}(\mU)=\begin{pmatrix}0&1&0&0\\
\dfrac{\gamma-3}{2}u^2+\dfrac{\gamma-1}{2}v^2&(3-\gamma)u&(1-\gamma)v&\gamma-1\\
-uv&v&u&0\\-\dfrac{\gamma uE}{\rho}+(\gamma-1)u(u^2+v^2)&H-(\gamma-1)u^2&(1-\gamma)uv&\gamma u\end{pmatrix},
\end{aligned}
$$
where, as in the 1-D case, $H=\frac{E+p}{\rho}$. We then introduce the matrices
$$
\begin{aligned}
\widehat A_{\jph,k}=\begin{pmatrix}0&1&0&0\\
\dfrac{\gamma-3}{2}\hat u^2+\dfrac{\gamma-1}{2} \hat v^2&(3-\gamma)\hat u&(1-\gamma)\hat v&\gamma-1\\
-\hat u\hat v&\hat v&\hat u&0\\-\dfrac{\gamma\hat u\hat E}{\hat\rho}+(\gamma-1)\hat u(\hat u^2+\hat v^2)&\hat H-(\gamma-1)\hat u^2&
(1-\gamma)\hat u\hat v&\gamma\hat u\end{pmatrix},
\end{aligned}
$$
where $\hat{(\cdot)}$ stands for the following averages:
\begin{equation*}
\begin{aligned}
&\hat\rho=\frac{\xbar\rho_{j,k}+\xbar\rho_{j+1,k}}{2},\quad\hat u=\frac{u_{j,k}+u_{j+1,k}}{2},\quad\hat v=\frac{v_{j,k}+v_{j+1,k}}{2},\quad
\hat p=\frac{p_{j,k}+p_{j+1,k}}{2},\\
&\hat H=\frac{\hat E+\hat p}{\hat\rho},\quad
\hat E=\frac{\hat p}{\gamma-1}+\frac{\hat\rho}{2}(\hat u^2+\hat v^2),
\end{aligned}
\end{equation*}
where
$$
u_{j,k}=\frac{(\xbar{\rho u})_{j,k}}{\xbar\rho_{j,k}},\quad v_{j,k}=\frac{(\xbar{\rho v})_{j,k}}{\xbar\rho_{j,k}},\quad\mbox{and}\quad
p_{j,k}=(\gamma-1)\left[\xbar E_{j,k}-\frac{\rho_{j,k}}{2}\left(u_{j,k}^2+v_{j,k}^2\right)\right].
$$
As in Appendix \ref{appa}, all of the $\hat{(\cdot)}$ quantities have to have a subscript index, that is,
$\hat{(\cdot)}=\hat{(\cdot)}_{\jph,k}$, which we omit for the sake of brevity for all of the quantities except for $\widehat A$. We then
compute the matrix $R_{\jph,k}$ composed of the right eigenvectors of $\widehat A_{\jph,k}$ and obtain
$$
R_{\jph,k}=\begin{pmatrix}1&1&0&1\\\hat u-\hat c&\hat u&0&\hat u+\hat c\\\hat v&\hat v&1&\hat v\\
\hat H-\hat u\hat c&\dfrac{\hat u^2+\hat v^2}{2}&\hat v&\hat H+\hat u\hat c\end{pmatrix}
$$
and
$$
R^{-1}_{\jph,k}=\dfrac{1}{\hat\phi}\begin{pmatrix}\dfrac{\hat u^2+\hat v^2}{2}+\dfrac{\hat u\hat\phi}{2\hat c}&
-\hat u-\dfrac{\hat\phi}{2\hat c}&-\hat v&1\\[1.2ex]2\hat\phi-2\hat H&2\hat u&2\hat v&-2\\-\hat v\hat\phi&0&\hat\phi&0\\
\dfrac{\hat u^2 +\hat v^2}{2}-\dfrac{\hat u \hat \phi }{2\hat c}&-\hat u+\dfrac{\hat\phi}{2\hat c}&-\hat v&1\end{pmatrix},
$$
where $\hat\phi=2\hat H-\hat u^2-\hat v^2$ and $\hat c=\sqrt{\gamma\hat p/\hat\rho}$.

We then introduce the local characteristic variables in the neighborhood of $(x,y)=(x_\jph,y_k)$:
\begin{equation*}
\bm\Gamma_{\ell,k}=R^{-1}_{\jph,k}\xbar\mU_{\ell,k},\quad\ell=j-1,\ldots,j+2.
\end{equation*}
Equipped with the values $\bm\Gamma_{j-1,k}$, $\bm\Gamma_{j,k}$, $\bm\Gamma_{j+1,k}$, and $\bm\Gamma_{j+2,k}$, we compute
\begin{equation*}
(\bm\Gamma_x)_{j,k}={\rm minmod}\left(2\,\frac{\bm\Gamma_{j,k}-\bm\Gamma_{j-1,k}}{\dx},\,\frac{\bm\Gamma_{j+1,k}-\bm\Gamma_{j-1,k}}{2\dx},\,
2\,\frac{\bm\Gamma_{j+1,k}-\bm\Gamma_{j,k}}{\dx}\right),
\end{equation*}
and
\begin{equation*}
(\bm\Gamma_x)_{j+1,\,k}={\rm minmod}\left(2\,\frac{\bm\Gamma_{j+1,k}-\bm\Gamma_{j,k}}{\dx},\,
\frac{\bm\Gamma_{j+2,k}-\bm\Gamma_{j,k}}{2\dx},\,2\,\frac{\bm\Gamma_{j+2,k}-\bm\Gamma_{j+1,k}}{\dx}\right),
\end{equation*}
where the minmod function, defined in \eref{2.7}, is applied in the component-wise manner. We then use these slopes to evaluate
$$
\bm\Gamma^{\rm E}_{j,k}=\,\bm\Gamma_{j,k}+\frac{\dx}{2}(\bm\Gamma_x)_{j,k}\quad\mbox{and}\quad
\bm\Gamma^{\rm W}_{j+1,k}=\bm\Gamma_{j+1,k}-\frac{\dx}{2}(\bm\Gamma_x)_{j+1,k},
$$
and finally obtain the corresponding point values of $\mU$ by
\begin{equation*}
\mU^{\rm E,W}_\jph=R_{\jph,k}\bm\Gamma^{\rm E,W}_{\jph,k}.
\end{equation*}

\bibliographystyle{siam}
\bibliography{ref}
\end{document}